\documentclass[12pt]{article}
\usepackage{amssymb}
\usepackage{amsfonts}
\usepackage{graphicx}
\usepackage{amsmath}
\title{Sharp regularity of linearization for $C^{1,1}$ hyperbolic diffeomorphisms
\thanks{Supported by NSFC and MOE research grants}
}

\author{
{\bf Wenmeng Zhang}\,$^{a,b}$, ~
{\bf Weinian Zhang}\,$^{a}$
\thanks{Corresponding to: matzwn@126.com or matwnzhang@yahoo.com.cn}
,
~
{\bf Witold Jarczyk}\,$^{c}$
\\
$^{a}${\small Yangtze Center of Mathematics and Department of Mathematics}\\
{\small Sichuan University, Chengdu, Sichuan 610064, P. R. China}
\\
$^{b}${\small College of Mathematics Science, Chongqing Normal University}\\
{\small Chongqing 400047, P. R. China}
\\
$^{c}${\small Faculty of Mathematics, Computer Science and Econometrics}
\\
{\small University of Zielona G\'ora, Szafrana 4a, 65-516 Zielona G\'ora, Poland}
}

\date{}

\begin{document}
\maketitle

\begin{abstract}
$C^1$ linearization is of special significance because it preserves
smooth dynamical behaviors and distinguishes qualitative properties
in characteristic directions. However, $C^1$ smoothness is not
enough to guarantee $C^1$ linearization. For $C^{1,1}$ hyperbolic diffeomorphisms
on Banach spaces $C^1$ linearization was proved under a gap condition together with a band
condition of the spectrum.
In this paper, the result of $C^1$ linearization in Banach spaces is strengthened to $C^{1,\beta}$ linearization
with a constant $\beta>0$ under a weaker band condition
by a decomposition with invariant foliations.
The weaker band condition allows the spectrum
to be a union of more than two but finitely many bands but restricts those bands to be bounded
by a number depending on
the supremum of contractive spectrum and the infimum of expansive spectrum.
Furthermore, we give an estimate for the exponent $\beta$
and prove that the estimate is sharp in the planar case.

\vspace{0.2cm}

\end{abstract}


\renewcommand{\theequation}{\thesection.\arabic{equation}}
\newtheorem{lm}{Lemma}
\newtheorem{thm}{Theorem}
\newtheorem{cor}{Corollary}
\newtheorem{re}{Remark}
\newtheorem{prop}{Proposition}
\allowdisplaybreaks

\vskip 0.4cm
\parskip 10pt


\section{Introduction}
\setcounter{equation}{0}

Let $X$ be a Banach space equipped with a norm $\|\cdot\|$ and let
$F:X\to X$ be a diffeomorphism fixing the origin $O$.
Linearizability of $F$ means that $F$ is conjugate to its linear
part $\Lambda:=DF(O)\in {\cal L}(X,X)$ near $O$, where ${\cal
L}(X,X)$ is the set of all bounded linear operators from $X$ into
$X$. Notice that the boundedness of $DF(O)$ comes from the
definition of (Fr\'{e}chet) differentiation (see \cite{ChangKC-book})
in Banach spaces.
The conjugacy is a solution of the functional equation
\begin{align}
\Phi\circ F=\Lambda\circ \Phi,
\label{schd-eqn}
\end{align}
where $\circ$ stands for the composition of mappings. Then
$F$ is said to be
$C^r$ linearizable (or to admit a $C^r$ linearization)
if equation (\ref{schd-eqn}) has a solution $\Phi$ which is a $C^r$ ($r\ge 0$) diffeomorphism.

A well-known result on linearization is the Hartman-Grobman's linearization theorem (cf. \cite{Hart64,pugh-AJM69}),
which tells that
every $C^1$ diffeomorphism can be $C^0$ linearized near a hyperbolic fixed point,
at which the spectrum of the linear part does not intersect the unit circle.
Efforts are also made to smooth linearization, i.e.,
the question:
Can a $C^k$ diffeomorphism be $C^r$ ($1\le r\le k$) linearized and how large can the order $r$ be?
The smooth linearization was initiated
by Sternberg in $\mathbb{R}^n$ in 1950's (cf. \cite{ster57, ster58}).
He proved that the number $r$ depends on both $k$ and
the non-resonance condition (cf. \cite{CLW-book94}), i.e., the eigenvalues $\lambda_1,...,\lambda_n$ of $\Lambda$
satisfy
\begin{align*}
\lambda_j\ne \lambda_1^{k_1}\cdots\lambda_n^{k_n}, ~~~~~ \forall
j=1,...,n,
\end{align*}
for any nonnegative integers $k_1,...,k_n$ such that $
\sum_{i=1}^{n}k_i\ge 2$.
Under this condition some improvements were made in e.g. \cite{Bronkopa-book94,Sell-AJM85} for $C^r$
linearization with $1\le r <\infty$, but all of those results involve a loss of smoothness,
i.e., $r<k$. The loss suggests an interesting work to find optimal results for smaller or minimal loss.

The smoothness of linearization is important to preserve more
dynamical properties in the procedure of linearization, such as the
smoothness of invariant manifolds and the characteristic directions
of the systems. In particular, $C^1$ linearization is of special
interest because it is useful to the problems of homoclinic
bifurcations (cf. \cite{AAIS-book94,DBo-JDE89}), stability of
topological mixing
for hyperbolic flows (cf. \cite{FMT-Ann07}), $C^1$ iterative roots
of mappings (cf. \cite{ZZJZ-JMAA12}), etc. In addition, $C^{1,\beta}$ ($\beta>0$) linearization
is also important, for example, to the studies of Lorenz attractors
(cf. \cite{HollMelb-07,Rob-SIAM-JMA92}). Here, for an
integer $r\ge 0$ and a real $\epsilon\in(0,1]$, $C^{r,\epsilon}$
denotes the class of all $C^r$ mappings $h$'s satisfying that
\begin{align}
\sup_{x\ne \tilde{x}}\frac{\|D^rh(x)-D^rh(\tilde{x})\|}{\|x-\tilde{x}\|^{\epsilon}}<\infty,
\label{Cnorm}
\end{align}
where $D^r$ denotes the $r$-th order differential operator.

In what follows we pay attention to $C^{1,1}$ mappings because
it was shown in \cite[p.139]{Kuczma1968} that there exists a $C^1$ contraction which cannot be $C^1$ linearized
even in one-dimensional case. In 1960, by investigating functional equations deduced from the
conjugacy, Hartman \cite{Hart60} proved that all $C^{1,1}$
contractions on $\mathbb{R}^n$ admit $C^{1,\beta}$ linearization
with $\beta>0$ depending on eigenvalues of $\Lambda$. Note that
expansions can be discussed by considering their inverses. In
2000's, using Hartman's idea, the result of $C^1$ linearization for
contractions was generalized to Banach spaces. For convenience,
let $\sigma(\Lambda)$ denote the spectrum of the linear part $\Lambda$
and suppose that
\begin{align}
\sigma(\Lambda)=\sigma_1\cup\,\cdots\,\cup\, \sigma_m,
\label{spec-decomp}
\end{align}
where $m\in\mathbb{N}$ and, for each $i=1,...,m$, $\sigma_i\subset \mathbb{C}$ and the numbers
$
\lambda_i^-:=\inf\{|z|:z\in \sigma_i\}
$
and
$
\lambda_i^+:=\sup\{|z|:z\in \sigma_i\}
$
satisfy
\begin{align}
0<\lambda_1^-\le\lambda_1^+<\cdots<\lambda_{m-1}^-\le\lambda_{m-1}^+<\lambda_m^-\le\lambda_m^+<1.
\label{spec-order}
\end{align}
As indicated in \cite{El-JMAA99,HPS-book77}, the $C^1$ smoothness of invariant manifolds
is guaranteed by {\it gaps} in the spectrum $\sigma(\Lambda)$, which can be formulated by the ratio $\lambda_i^-/\lambda_{i-1}^+$
for each $i\in\{2,...,m\}$. Those manifolds are crucial in the proof of linearization
as shown in \cite{El-JFA01,R-S-JDE04} because they are employed to construct transformations
to simplify the mapping $F$.
It is also required that the {\it band width}
of each spectral subset $\sigma_i$, which can be formulated by the ratio $\lambda_i^+/\lambda_{i}^-$ for each $i\in\{1,...,m\}$,
is so small that the well-known Banach's contraction principle can be applied to
each functional equation deduced from the conjugacy.
More precisely,
it was proved in \cite{El-JFA01} that $C^{1,1}$ contractions admit
$C^{1,\beta}$ linearization with a constant $\beta>0$ if there are real
$s_i\in(0,1)$, $i=1,...,m$, such that
\begin{align}
\lambda_i^+/\lambda_i^-<(\lambda_m^+)^{-\delta_i} ~~~ {\rm and} ~~~
(\lambda_i^+)^{s_i}/\lambda_i^-<(\lambda_m^+)^{-(1+\delta_{i-1})},
~~~\forall i=1,...,m,
\label{bond-cond}
\end{align}
where $\delta_m=1$ and $\delta_j=\delta_{j+1}(1-s_{j+1})$ for all
$j=0,...,m-1$.
Later,
in \cite{R-S-JDE04} condition (\ref{bond-cond}) was replaced with the weaker one
\begin{align}
\lambda_i^+/\lambda_i^-<(\lambda_m^+)^{-1}, ~~~~~ \forall i=1,...,m,
\label{NR-conB}
\end{align}
i.e., band widths of all spectral subsets are uniformly less than
the reciprocal of $\lambda_m^+$, the supremum  of contractive
spectrum $\sigma(\Lambda)$, but the obtained result is $C^{1}$
linearization other than $C^{1,\beta}$ linearization.
Recently, for $C^1$ linearization of planar contractions, in
\cite{WZWZ-JFA01} the requirement of smoothness of $F$ was lowered from $C^{1,1}$ to $C^{1,\alpha}$
associated with the condition $\alpha> 1-\log |\lambda_2|/\log |\lambda_1|$
and the linearization is proved to be $C^{1,\beta}$ for some $\beta>0$.

The case that both contraction and expansion are included is more
complicated. Unlike the works \cite{El-JFA01,Hart60,R-S-JDE04} for
contractions, divergence factor in the present case cannot be avoid
no matter whether we consider $F$ or $F^{-1}$, which appears as an
obstacle to the applications of Banach's contraction principle. In
2004, extending Belitskii's result (\cite{Beli-FAA73,Beli-RMS78}) in
$\mathbb{R}^n$, which gives $C^1$ linearization under non-resonance
conditions, Rodrigues and Sol${\rm \grave{a}}$-Morales
(\cite{R-S-JDDE04}) obtained a result of $C^1$ linearization in
Banach spaces under a gap condition together with a band condition
of the spectrum. Concretely, with the spectral splitting
$\sigma(\Lambda)=\sigma_- \cup \sigma_+$, where
\begin{align}
\{|\lambda|\in\mathbb{R}:\lambda\in\sigma_-\}\subset
(\lambda_s^-,\lambda_s^+), ~~~~~~
\{|\lambda|\in\mathbb{R}:\lambda\in\sigma_+\}\subset
(\lambda_u^-,\lambda_u^+)
\label{tj}
\end{align}
and $0<\lambda_s^-<\lambda_s^+<1<\lambda_u^-<\lambda_u^+$,
it is proved in \cite{R-S-JDDE04} that a $C^{1,1}$ diffeomorphism on a Banach space can
be $C^1$ linearized if
$
\lambda_s^+\lambda_u^+<\lambda_u^-$, $\lambda_s^+<\lambda_u^-\lambda_s^-$,
$(\lambda_s^+)^2<\lambda_s^-$ and $\lambda_u^+<(\lambda_u^-)^2,
$
or equivalently if
\begin{align}
\left\{
\begin{array}{ll}
\lambda_u^-/\lambda_s^+>\max\{\lambda_u^+,(\lambda_s^-)^{-1}\},
\\
\vspace{-0.3cm}
\\
\lambda_s^+/\lambda_s^-<(\lambda_s^+)^{-1}, ~~~ \lambda_u^+/\lambda_u^-<\lambda_u^-.
\end{array}
\right.
\label{RS-NR}
\end{align}
The first inequality of (\ref{RS-NR}) indicates that the size of the spectral gap cannot
be too small and the next two inequalities in (\ref{RS-NR}) show that the band widths of
spectral subsets $\sigma_-$ and $\sigma_+$ are small enough. We note
that it is hard to strengthen the above-mentioned result from $C^1$
linearization to $C^{1,\beta}$ linearization with
$\beta>0$. In fact, in \cite{R-S-JDDE04} equation (\ref{schd-eqn})
was decomposed into three functional equations,
two of which depend on only the first variable and only the second variable separately
and the other of which is equal to 0 on both axes.
In order to solve the third equation, a special norm
\begin{align*}
\|h\|_{\alpha}:=\max\Big\{\sup_{x_2\ne 0}\frac{|\partial_{x_1}h(x_1,x_2)|}{|x_2|^\alpha},~
\sup_{x_1\ne 0}\frac{|\partial_{x_2}h(x_1,x_2)|}{|x_1|^\alpha}\Big\},
\end{align*}
where $\alpha\in(0,1]$, is employed so as to obtain the contraction constant
\begin{align*}
\iota:=\min\{\lambda_s^+(\lambda_u^+)^\alpha/\lambda_u^-,~~\lambda_u^+(\lambda_s^+)^\alpha/\lambda_u^-\}<1.
\end{align*}
The above inequality holds when $\alpha$ is sufficiently close to
$1$ due to the first inequality of (\ref{RS-NR}). So, the special
norm $\|\cdot\|_\alpha$ cannot be replaced with
$\|\cdot\|_{C^{1,\beta}}$, the norm defined by (\ref{Cnorm}) for
$r=1$ and $\epsilon=\beta$; otherwise the constant $\iota$
becomes $(\lambda_u^+)^{1+\beta}/\lambda_u^-$, which is always $>1$
because $1<\lambda_u^-<\lambda_u^+$ and therefore the Contraction
Principle is not applicable.

In this paper we study $C^{1,\beta}$ linearization of $C^{1,1}$
hyperbolic diffeomorphisms on Banach spaces. In order to overcome
the difficulties of the divergence factor, our strategy consists of
the following three steps. First, we decompose $F$ along invariant
foliations, which can be constructed by solving the Lyapunov-Perron
equation, and reduce the problem to the linearization of a
contraction and the linearization of an expansion separately. Unlike
\cite{KirPal90, Lu-JDE91, Tan-JDE00}, where $C^0$ linearization and
$C^{0,\alpha}$ linearization are discussed, we need to consider a
system of functional equations, one of which comes from derivatives
of the Lyapunov-Perron equation, for $C^{1,\beta}$ solutions of the
Lyapunov-Perron equation. Moroever, a special technique  (see the
proof of our Lemma \ref{lm-hest}) will be used in estimating the
H\"{o}lder exponent $\beta$. Second, in order to prove $C^1$
linearization of a contraction, we use a sequence of mappings where
iterates of $F$ are involved to approach the solution of equation
(\ref{schd-eqn}). This result itself improves the corresponding ones
given in \cite{El-JFA01,R-S-JDE04} due to a weaker condition on
$\sigma(\Lambda)$ and a larger bound of $\beta$ (see Remark
\ref{re-comNRcc} below). Third, we give a solution of equation
(\ref{schd-eqn}) by composing the two transformations obtained in
the previous two steps separately and conclude our result
Theorem~\ref{thm-linB} of linearization in the general hyperbolic
case. This result extends Rodrigues and Sol${\rm
\grave{a}}$-Morales' result given in \cite{R-S-JDDE04} by weakening
condition (\ref{RS-NR}) (see Remark \ref{re-comNR} below) and
strengthening $C^1$ smoothness of linearization to $C^{1,\beta}$.
Moreover, an estimate of $\beta$ is given and proved to be sharp in
the planar case.

Our paper is
organized as follows: After introducing some basic properties of $F$,
useful notations and the fiber contraction theorem in Section 2,
we investigate solutions of the Lyapunov-Perron equation,
proving the existence of $C^1$ solutions in Section 3 and
showing  in Section 4 that those solutions are actually $C^{1,\beta}$ with some real $\beta>0$.
Then,
we study partial linearization for contractions in Section 5 and leave proofs of some lemmas to
Section 6.
In Section 7 we employ the obtained result of partial linearization to prove $C^{1,\beta}$ linearization for contractions.
Finally, combining the above obtained results, we give our theorem on $C^{1,\beta}$ linearization
for hyperbolic mappings in Section 8.
In Section 9 we show that the estimate of $\beta$ is sharp in the planar case.


\section{Preliminaries}
\setcounter{equation}{0}

In the following, we discuss on a Banach space $(X, \|\cdot\|)$
which has bump functions (smooth functions on $X$ each of which has a bounded support, i.e.,
the function is identical to $1$ within the support but to $0$ outside a neighborhood of the support).
Known from \cite{FryMc-EM2, LeWh-PAMS79},
every Hilbert space has bump functions but some Banach spaces (e.g., $C^0[0,1]$) do not have a bump function.
Throughout this paper, we always let $K,L,M$ and $K_i,L_i,M_i$, $\forall i\in \mathbb{N}$, be positive constants.

Let $\varrho$ be a bump function on $X$ such that
\begin{align*}
\varrho(x)=1, ~~~~~\forall x\in U, ~~~~~{\rm and}~~~~~
\varrho(x)=0, ~~~~~ \forall x\in X\backslash V,
\end{align*}
where $U$ and $V$, $U\subsetneq V \subset X$, are neighborhoods of
the origin $O$. Without loss of generality, we may assume that $U$
and $V$ are sufficiently small; otherwise, we can find a smooth
transformation (e.g., dilation) to change $U$ and $V$ to be sufficiently
small ones and send $\varrho$ to another bump function. Multiplying
the higher order terms of $F$ by $\varrho$, we obtain a modified
diffeomorphism such that
\begin{align}
\|DF(x)-\Lambda\|\le \eta, \quad \forall x\in X,
\label{FFbumpR2}
\end{align}
where $\eta>0$ is a sufficiently small constant provided $V$ is
small enough. This modification does not affect our results at all
because we are only interested in local properties of $F$.

Next, assume that spectrum $\sigma(\Lambda)$ is an union of two sets $\sigma_-$ and $\sigma_+$ such that (\ref{tj}) holds.
By the Spectral Decomposition Theorem
(see, e.g., \cite[p.9]{GGKbook-90}) one can further assume that the space $X$ has a direct sum decomposition
$
X=X_-\oplus X_+
$
with $\Lambda$-invariant subspaces $X_-$ and $X_+$, that is,
\begin{align}
\Lambda={\rm diag}(\Lambda_-, \Lambda_+),
\label{def-FB}
\end{align}
where $\Lambda_-\in {\cal L}(X_-, X_-)$ and $\Lambda_+\in {\cal L}(X_+, X_+)$
such that
$
\sigma(\Lambda_-)=\sigma_-
$
and
$
\sigma(\Lambda_+)=\sigma_+
$
respectively. Moreover, we let
\begin{align}
\sigma_-=\bigcup_{i=1}^{d}~\sigma_i
\quad{\rm and}\quad
\sigma_+=\bigcup_{i=d+1}^{m}\sigma_i,
\label{spec-cond}
\end{align}
where the integer $m\ge 2$, $d\in\{1,...,m-1\}$, $\sigma_i\subset \mathbb{C}$ and the numbers
$\lambda_i^-:=\inf\{|z|:z\in \sigma_i\},~
\lambda_i^+:=\sup\{|z|:z\in \sigma_i\},~\forall i=1,...,m,$
satisfy
\begin{align}
0<\lambda_1^-\le \lambda_1^+<\cdots<\lambda_d^-\le \lambda_d^+<
1<\lambda_{d+1}^-\le
\lambda_{d+1}^+<\cdots<\lambda_m^-\le\lambda_m^+.
\label{spec-distB}
\end{align}
Then (\ref{tj}) yields
$\lambda_s^-<\lambda_1^-$, $\lambda_d^+<\lambda_s^+$,
$\lambda_u^-<\lambda_{d+1}^-$, $\lambda_m^+<\lambda_u^+$.
Without loss of generality, we put
\begin{align*}
\lambda_s^-=\lambda_1^- -\delta,
\quad
\lambda_s^+=\lambda_d^+ +\delta,
\quad
\lambda_u^-=\lambda_{d+1}^- -\delta,
\quad
\lambda_u^+=\lambda_m^+ +\delta
\end{align*}
for a sufficiently small number $\delta>0$. Obviously, we still have $0<\lambda_s^-<\lambda_s^+<1<\lambda_u^-<\lambda_u^+$.
By \cite[Theorem 5]{R-S-JDE04}, for any number
$\delta>0$, one can choose appropriate equivalent norms
in $X_-$ and $X_+$ separately such that
\begin{align}
\|\Lambda_-\|< \lambda_s^+, \quad \|\Lambda_-^{-1}\|<1/\lambda_s^-,  \quad
\|\Lambda_+\|< \lambda_u^+, \quad \|\Lambda_+^{-1}\|< 1/\lambda_u^-.
\label{mu-lam}
\end{align}

The following notations are useful.
Let $C^0_b(\Omega,Z_2)$ consist of all $C^0$ mappings $h$ from $\Omega$, an open subset of a Banach space $(Z_1,\|\cdot\|_1)$, into
another Banach space $(Z_2,\|\cdot\|_2)$ such that $\sup_{z\in \Omega}\|h(z)\|_2<\infty$.
It is a Banach space equipped with the norm $\|\cdot\|_{C^0_b(\Omega,Z_2)}$ defined by
\begin{align*}
\|h\|_{C^0_b(\Omega,Z_2)}:=\sup_{z\in \Omega}\|h(z)\|_2, ~~~~~ \forall h\in C^0_b(\Omega,Z_2).
\end{align*}
For a positive number $\gamma\in\mathbb{R}$, let $S_\gamma(\Omega,Z_2)$ consist of all sequences $u:=(u_n)_{n\ge 0}$, where $u_n\in
C^0_b(\Omega,Z_2)$ for all integers $n\ge 0$, such that
$
\sup_{n\ge 0}\{\gamma^{-n}\|u_n\|_{C^0_b(\Omega,Z_2)}\}<\infty.
$
Then, $S_\gamma(\Omega,Z_2)$ is a Banach space equipped with the norm $\|\cdot\|_{S_{\gamma}(\Omega,Z_2)}$
defined by
\begin{align*}
\|u\|_{S_{\gamma}(\Omega,Z_2)}:=\sup_{n\ge 0}\{\gamma^{-n}\|u_n\|_{C^0_b(\Omega,Z_2)}\},
\qquad \forall u\in S_{\gamma}(\Omega,Z_2).
\end{align*}
In fact,
suppose that $(u^{(k)})_{k\in\mathbb{N}}$, where $u^{(k)}:=(u^{(k)}_n)_{n\ge 0} \in S_{\gamma}(\Omega,Z_2)$,
is a Cauchy sequence in $S_{\gamma}(\Omega,Z_2)$. For any given $\epsilon>0$, there is an integer $K_0>0$
such that
\begin{align}
\|u^{(k_1)}-u^{(k_2)}\|_{S_{\gamma}(\Omega,Z_2)}=\sup_{n\ge 0}\{\gamma^{-n}\|u^{(k_1)}_n-u^{(k_2)}_n\|_{C^0_b(\Omega,Z_2)}\}<\epsilon,
\label{xnyym}
\end{align}
provided that $k_1,k_2>K_0$. This implies that $(u^{(k)}_n)_{k\in\mathbb{N}}$ is also a Cauchy sequence in $C^0_b(\Omega,Z_2)$
for each $n\ge 0$ and consequently there exists an $u^*_n\in C^0_b(\Omega,Z_2)$ such that
$
\lim_{k\to \infty}\|u^{(k)}_n-u^*_n\|_{C^0_b(\Omega,Z_2)}=0.
$
Thus, putting $u^*:=(u^*_n)_{n\ge 0}$ and letting $k_2\to \infty$ in (\ref{xnyym}), one easily verifies that
$u^*\in {S_{\gamma}(\Omega,Z_2)}$
and
$\lim_{k\to \infty}\|u^{(k)}-u^*\|_{S_{\gamma}(\Omega,Z_2)}=0$.
This shows that $S_\gamma(\Omega,Z_2)$ is a Banach space equipped with the norm $\|\cdot\|_{S_{\gamma}(\Omega,Z_2)}$.

The following lemma is the Fiber Contraction Theorem, which will be used in this paper and
can be found in \cite{HirsPugh-70} (see also \cite[p.111]{Chicone-book99}).

\begin{lm}
Let ${\cal T}:Z_1\to Z_1$ and ${\cal S}:Z_1\times Z_2\to Z_2$ be mappings and let ${\cal Q}:Z_1\times Z_2\to Z_1\times Z_2$
be defined by
\begin{align*}
{\cal Q}(z_1,z_2):=({\cal T}(z_1),{\cal S}(z_1,z_2)), ~~~~~ \forall (z_1,z_2)\in Z_1\times Z_2.
\end{align*}
Suppose that ${\cal T}$ is a contraction and that ${\cal Q}$ is a continuous fiber contraction, i.e.,
${\cal Q}$ is continuous such that
\begin{align*}
\|{\cal Q}(z_1,z_2)-{\cal Q}(z_1,\tilde{z}_2)\|\le \theta\|z_2-\tilde{z}_2\|_2,
~~~~~ \forall z_1\in Z_1, ~ \forall z_2,\tilde{z}_2\in Z_2,
\end{align*}
with $\theta\in (0,1)$. Then $(z_1^*,z_2^*)\in Z_1\times Z_2$ is an attracting
fixed point of ${\cal Q}$, i.e.,
\begin{align*}
\lim_{n\to \infty}{\cal Q}^n(z_1,z_2)=(z_1^*,z_2^*),
~~~~~ \forall (z_1,z_2)\in Z_1\times Z_2,
\end{align*}
where $z_1^*\in Z_1$ is the fixed point of ${\cal T}$ and $z_2^*\in Z_2$ is the fixed point of ${\cal S}(z_1^*,\cdot)$.
\label{lm-fbcont}
\end{lm}



\section{$C^1$ solution of Lyapunov-Perron equation}
\setcounter{equation}{0}

Let $f=F-\Lambda$ be the nonlinear term of $F$ and let $\pi_-$ and
$\pi_+$ be projections onto $X_-$ and $X_+$ respectively.
As indicated in the Introduction, we need to decompose $F$ along invariant foliations.
According to \cite[p.117]{HPS-book77} and \cite{ChHaTan-JDE97},
a {\it stable invariant foliation} of $X$ for $F$ is a disjoint decomposition of $X$ into injectively immersed
connected submanifolds called {\it leaves}, which can be formulated by
\begin{align}
{\cal M}_s(x)=\big\{\big(\,y_-,h_s(x,y_-)\big): y_-\in X_-\big\}, ~~~~~ \forall x\in X,
\label{def-Ms}
\end{align}
for some mappings $h_s:X\times X_-\to X_+$ with $h_s(x,\pi_-x)=\pi_+x$, such that
\begin{align*}
F({\cal M}_s(x))\subset {\cal M}_s(F(x))
\end{align*}
(see Fig.1). A definition of unstable invariant foliation $\{{\cal M}_u(x)\}_{x\in X}$ of $X$ for $F$ can be stated analogously.


~~~~~~~~~~~~~~~~~~~~~~~~~~~~~~~~~~~~
\begin{picture}(390,180)(0,0)
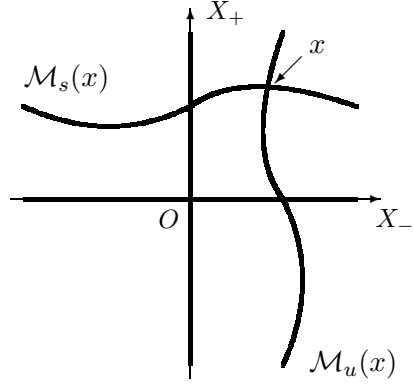
\put(-15,5){\footnotesize Fig.1. Leaves of stable and unstable invariant foliations}
\put(27,95){\vector(1,0){140}}                                      
\put(165,86){\makebox(2,1)[l]{\footnotesize $X_-$}}
\put(95,27){\vector(0,1){140}}                                      
\put(101,165){\makebox(2,1)[l]{\footnotesize $X_+$}}
\put(83,87){\makebox(2,1)[l]{{\footnotesize $O$}}}                  

\put(140,152){\makebox(2,1)[l]{\small $x$}}
\put(137,149){\vector(-1,-1){10}}


{\linethickness{0.5mm}
\put(32,95){\line(1,0){126}}                                        
}


{\linethickness{0.5mm}
\qbezier(32,130)(65,115)(95,130)                                    
\qbezier(95,130)(115,145)(158,130)
}

\put(32,140){\makebox(2,1)[l]{\small ${\cal M}_s(x)$}}

{\linethickness{0.5mm}
\put(95,32){\line(0,1){126}}                                        
}

{\linethickness{0.5mm}
\qbezier(130,32)(145,65)(130,95)                                       
\qbezier(130,95)(115,115)(130,158)
}

\put(140,32){\makebox(2,1)[l]{\small ${\cal M}_u(x)$}}

\end{picture}


In order to construct invariant foliations, we need to study the Lyapunov-Perron equation (cf. \cite{ChHaTan-JDE97})
\begin{align}
q_n(x,y_-)=~&\Lambda_-^n(y_--\pi_- x)+\sum_{k=0}^{n-1}\Lambda_-^{n-k-1}\{\pi_- f(q_k(x,y_-)+F^k(x))-\pi_- f(F^k(x))\}
\nonumber\\
&-\sum_{k=n}^{\infty}\Lambda_+^{n-k-1}\{\pi_+ f(q_k(x,y_-)+F^k(x))-\pi_+ f(F^k(x))\}, ~~~ \forall n\ge 0,
\label{eqns-foli}
\end{align}
where $q_n:X\times X_-\to X$ is unknown for every integer $n\ge 0$.
In fact, assume that $(q_n(x,y_-))_{n\ge 0}$ is a sequence such that
$\sup_{n\ge 0}\{\gamma^{-n}\|q_n(x,y_-)\|\}< \infty$ holds with
$\gamma\in (\lambda_s^+,\lambda_u^-)$ for every point $(x,y_-)\in
X\times X_-$. Then one verifies straightforward that
$(q_n(x,y_-))_{n\ge 0}$ is a solution of equation (\ref{eqns-foli})
if and only if
\begin{align}
F^n(x+q_0(x,y_-))-F^n(x)=q_n(x,y_-),~~~~~ \forall n\ge 0.
\label{n--0}
\end{align}
Therefore, if the above-mentioned solution $(q_n(x,y_-))_{n\ge 0}$
is unique then
\begin{align*}
{\cal M}_s(x):=\{x+q_0(x,y_-):y_-\in X_-\} =\{\tilde{x}\in X:
\sup_{n\ge 0}\{\gamma^{-n}\|F^n(\tilde{x})-F^n(x)\|\}< \infty\}
\end{align*}
for every $x\in X$, which is obviously invariant under $F$.
Moreover, noting that $\pi_-(x+q_0(x,y_-))=y_-$ and putting
$h_s(x,y_-)=\pi_+(x+q_0(x,y_-))$ we get (\ref{def-Ms}). Equation
(\ref{eqns-foli}) was discussed in \cite{ChHaTan-JDE97} for $C^0$
solution $(q_n)_{n\ge 0}$ such that each $q_n(x,\cdot): X_-\to X$ is
$C^1$, i.e., each leaf of the invariant foliation is $C^1$. For our
purpose of $C^1$ linearization, in this paper we further need those
$C^1$ leaves to be jointed smoothly. Therefore, in this section we
find $C^1$ solution $(q_n)_{n\ge 0}$ of equation (\ref{eqns-foli}),
i.e., each $q_n: X\times X_-\to X$ is $C^1$.

\begin{lm}
Let $\Omega\subset X\times X_-$ be a small neighborhood of the origin $\tilde{O}\in X\times X_-$.
Suppose that $F: X\to X$ is a $C^{1,1}$ diffeomorphism satisfying {\rm (\ref{FFbumpR2})}
and that {\rm (\ref{mu-lam})} holds
with
\begin{align}
\lambda_s^+\lambda_u^+<\lambda_u^- ~~~~~ {\rm(\,or} ~
\lambda_d^+\lambda_m^+<\lambda_{d+1}^- {\rm)}.
\label{NR22}
\end{align}
Then the Lyapunov-Perron equation {\rm (\ref{eqns-foli})} has a
unique solution $(q_n)_{n\ge 0}$ such that every $q_n:\Omega\to X$
for $n\ge 0$ is of class $C^1$ and
\begin{align*}
(q_n)_{n\ge 0}\in S_{\gamma_1}(\Omega,X), \qquad
(Dq_n)_{n\ge 0}\in S_{\gamma_2}(\Omega,{\cal L}(X\times X_-,X)),
\end{align*}
where the positive numbers $\gamma_1,\gamma_2$ satisfy
\begin{align}
\lambda_s^+<\gamma_1<1<\gamma_2<\lambda_u^- \quad{\rm and}\quad \gamma_1\lambda_u^+<\gamma_2.
\label{gama12}
\end{align}
\label{lm-LPeqn}
\end{lm}

\vspace{-0.8cm}

Note that the choice of $\gamma_1$ and $\gamma_2$ in the above lemma is reasonable due to inequality (\ref{NR22}).

{\bf Proof}. In the proof, we let
\begin{align*}
\|x\|:=\|x_-\|+\|x_+\|, ~~~~~~~ \|(x,y_-)\|:=\|x\|+\|y_-\|,
\end{align*}
for $x=x_-+x_+\in X$ and $y_-\in X_-$. Let
$E_1:=S_{\gamma_1}(\Omega,X)$ and $E_2:=S_{\gamma_2}(\Omega,{\cal L}(X\times X_-,X))$
for short.
As defined in Section 2, $E_1$ and $E_2$ are both Banach spaces equipped the corresponding norms,
denoted by $\|\cdot\|_{E_1}$ and $\|\cdot\|_{E_2}$ respectively.
Define operators ${\cal T}: E_1\to E_1$ and ${\cal S}: E_1\times E_2\to E_2$ by
\begin{align}
({\cal T} v)_n(x,y_-):=\,&\Lambda_-^n(y_--\pi_- x)
\nonumber\\
&+\sum_{k=0}^{n-1}\Lambda_-^{n-k-1}\{\pi_- f(v_k(x,y_-)+F^k(x))-\pi_-
f(F^k(x))\}
\nonumber\\
&-\sum_{k=n}^{\infty}\Lambda_+^{n-k-1}\{\pi_+ f(v_k(x,y_-)+F^k(x))-\pi_+
f(F^k(x))\}
\label{def-T}
\end{align}
and
\begin{align}
{\cal S}(v,w)_n(x,y_-)
:=\,&\Big({\rm
diag}(0,-\Lambda_-^n),\Lambda_-^n\Big)
\nonumber\\
&+\sum_{k=0}^{n-1}\Lambda_-^{n-k-1}\{D(\pi_-
f)(v_k(x,y_-)+F^k(x))(w_k(x,y_-)+DF^k(x))
\nonumber\\
&\hspace{2.8cm}-D(\pi_- f)(F^k(x))DF^k(x)\}
\nonumber\\
&-\sum_{k=n}^{\infty}\Lambda_+^{n-k-1}\{D(\pi_+
f)(v_k(x,y_-)+F^k(x))(w_k(x,y_-)+DF^k(x))
\nonumber\\
&\hspace{2.8cm}-D(\pi_+ f)(F^k(x))DF^k(x)\}
\label{def-S}
\end{align}
respectively for all $v:=(v_n)_{n\ge 0}\in E_1$ and all $w:=(w_n)_{n\ge 0}\in
E_2$. Note that in the above formula (\ref{def-S}) we have
\begin{align*}
\Big({\rm diag}(0,-\Lambda_-^n),~\Lambda_-^n\Big)\in {\cal L}(X\times X_-,X).
\end{align*}
Clearly, a fixed point of ${\cal T}$ is a solution of equation (\ref{eqns-foli}).
However, just applying the Contraction Principle to ${\cal T}$ is not enough because finding $C^1$ solutions of equation (\ref{eqns-foli})
requires us to deal with (\ref{eqns-foli}) itself together with its derivative. Thus we need to apply Lemma \ref{lm-fbcont}
(Fiber Contraction Theorem)
to both ${\cal T}$ and ${\cal S}$ in what follows.
For this purpose, we claim the following:
\begin{description}
\item (A1)
The operators ${\cal T}$ and ${\cal S}$ are well defined.
\item (A2)
The operator ${\cal Q}: E_1\times E_2\to E_1\times E_2$ defined by
\begin{align}
{\cal Q}(v,w):=({\cal T}v, {\cal S}(v,w)), ~~~~~ \forall (v,w)\in
E_1\times E_2,
\label{def-Q}
\end{align}
has an attracting fixed point $(v_*,w_*)\in E_1\times E_2$, i.e.,
\begin{align}
\lim_{n\to \infty}{\cal Q}^n(v,w)=(v_*,w_*),
~~~~~ \forall (v,w)\in E_1\times E_2,
\label{limQ}
\end{align}
where $v_*\in E_1$ is the fixed point of ${\cal T}$ and
$w_*\in E_2$ is the fixed point of ${\cal S}(v_*,\cdot)$.
\label{lm-fc}
\end{description}

Assertion (A1) holds because (\ref{FFbumpR2}) and (\ref{gama12}) yield
\begin{align*}
&\gamma_1^{-n}\|({\cal T} v)_n(x,y_-)\|
\\
&\le \Big(\frac{\lambda_s^+}{\gamma_1}\Big)^{
n}\|(x,y_-)\|+\gamma_1^{-n}\sum_{k=0}^{n-1}\|\Lambda_-^{n-k-1}\|\,\|\pi_-
f(v_k(x,y_-)+F^k(x))-\pi_- f(F^k(x))\|
   \\
   &~~~~+\gamma_1^{-n}\sum_{k=n}^{\infty}\|\Lambda_+^{n-k-1}\|\,\|\pi_+ f(v_k(x,y_-)+F^k(x))-\pi_+ f(F^k(x))\|
\\
&\le
1+\gamma_1^{-1}\sum_{k=0}^{n-1}\Big(\frac{\lambda_s^+}{\gamma_1}\Big)^{
n-k-1}
    \sup_{\xi\in X}\|D(\pi_- f)(\xi)\|\,\gamma_1^{-k}\|v_k(x,y_-)\|
   \\
   &~~~~+\gamma_1^{-1}\sum_{k=n}^{\infty}\Big(\frac{\gamma_1}{\lambda_u^-}\Big)^{k+1-n}
    \sup_{\xi\in X}\|D(\pi_+ f)(\xi)\|\,\gamma_1^{-k}\|v_k(x,y_-)\|
\\
&\le 1+K\|v\|_{E_1},
\end{align*}
implying that ${\cal T} v:=(({\cal T} v)_n)_{n\ge 0}\in E_1$, and
\begin{align*}
&\gamma_2^{-n}\|{\cal S}(v, w)_n(x,y_-)\|
\\
&\le\Big(\frac{\lambda_s^+}{\gamma_2}\Big)^n+\gamma_2^{-1}\sum_{k=0}^{n-1}\|\Lambda_-^{n-k-1}\|\Big(\|\{D(\pi_-
       f)(v_k(x,y_-)+F^k(x))\|\,\|w_k(x,y_-)\|
   \\
   &\hspace{0.5cm}~~~~+\|D(\pi_-f)(v_k(x,y_-)+F^k(x))-D(\pi_-f)(F^k(x))\|\,\|DF^k(x)\|\Big)
   \\
   &\hspace{0.5cm}+\gamma_2^{-1}\sum_{k=n}^{\infty}\|\Lambda_+^{n-k-1}\|\Big(\|\{D(\pi_+
       f)(v_k(x,y_-)+F^k(x))\|\,\|w_k(x,y_-)\|
   \\
   &\hspace{0.5cm}~~~~+\|D(\pi_+f)(v_k(x,y_-)+F^k(x))-D(\pi_+f)(F^k(x))\|\,\|DF^k(x)\|\Big)
\\
&\le\Big(\frac{\lambda_s^+}{\gamma_2}\Big)^n+\gamma_2^{-1}\sum_{k=0}^{n-1}\Big(\frac{\lambda_s^+}{\gamma_2}\Big)^{n-k-1}\Big(\sup_{\xi\in
X}
           \|D(\pi_- f)(\xi)\|\gamma_2^{-k}\|w_k(x,y_-)\|
   \\
   &\hspace{0.5cm}~~~~+L\gamma_2^{-k}\|v_k(x,y_-)\|\,\|DF^k(x)\|\Big)
   \\
   &\hspace{0.5cm}+\gamma_2^{-1}\sum_{k=n}^{\infty}\Big(\frac{\gamma_2}{\lambda_u^-}\Big)^{k+1-n}\Big(\sup_{\xi\in X}
           \|D(\pi_+ f)(\xi)\|\gamma_2^{-k}\|w_k(x,y_-)\|
   \\
   &\hspace{0.5cm}~~~~+L\gamma_2^{-k}\|v_k(x,y_-)\|\,\|DF^k(x)\|\Big)
\\
&\le1+\gamma_2^{-1}\sum_{k=0}^{n-1}\Big(\frac{\lambda_s^+}{\gamma_2}\Big)^{\,n-k-1}\Big(\eta\|w\|_{E_2}+
       L\Big(\frac{\gamma_1(\lambda_u^++\eta)}{\gamma_2}\Big)^{\!k\,}\|v\|_{E_1}\Big)
   \\
   &~~~~+\gamma_2^{-1}\sum_{k=n}^{\infty}\Big(\frac{\gamma_2}{\lambda_u^-}\Big)^{\,k+1-n}\Big(\eta\|w\|_{E_2}+
       L\Big(\frac{\gamma_1(\lambda_u^++\eta)}{\gamma_2}\Big)^{\!k\,}\|v\|_{E_1}\Big)
\\
&\le 1+K_1\|v\|_{E_1}+K_2\|w\|_{E_2},
\end{align*}
implying that ${\cal S}(v, w):=({\cal S}(v, w)_n)_{n\ge 0}\in E_2$.
Here we note
that $\|(x,y_-)\|<1$ since $\Omega$ is small and
that $\gamma_1(\lambda_u^++\eta)/\gamma_2<1$ by the second inequality of
(\ref{gama12}) since $\eta$ is small enough as mentioned just blow
(\ref{FFbumpR2}). In the above formula we understand that
$
\|DF^k(x)\|\le\|DF(F^{k-1}(x))\|\cdots\|DF(x)\|\le
(\lambda_u^++\eta)^k.
$

Assertion (A2) is proved by applying Lemma \ref{lm-fbcont} (Fiber Contraction Theorem).
First of all, ${\cal T}:E_1\to E_1$ is a contraction
because
\begin{align*}
&\gamma_1^{-n}\|({\cal T} v)_n(x,y_-)-({\cal T}
\tilde{v})_n(x,y_-)\|
\nonumber\\
&\le\gamma_1^{-n}\sum_{k=0}^{n-1}\|\Lambda_-^{n-k-1}\|\,
     \|\pi_- f(v_k(x,y_-)+F^k(x))-\pi_- f(\tilde{v}_k(x,y_-)+F^k(x))\|
\nonumber\\
&~~~~+\gamma_1^{-n}\sum_{k=n}^{\infty}\|\Lambda_+^{n-k-1}\|\,
     \|\pi_+ f(v_k(x,y_-)+F^k(x))-\pi_+ f(\tilde{v}_k(x,y_-)+F^k(x))\|
\nonumber\\
&\le\gamma_1^{-1}\sum_{k=0}^{n-1}\Big(\frac{\lambda_s^+}{\gamma_1}\Big)^{n-k-1}
     \sup_{\xi\in X}\|D(\pi_- f)(\xi)\|\gamma_1^{-k}\|v_k(x,y_-)-\tilde{v}_k(x,y_-)\|
\nonumber\\
&~~~~+\gamma_1^{-1}\sum_{k=n}^{\infty}\Big(\frac{\gamma_1}{\lambda_u^-}\Big)^{k+1-n}
     \sup_{\xi\in X}\|D(\pi_+ f)(\xi)\|\gamma_1^{-k}\|v_k(x,y_-)-\tilde{v}_k(x,y_-)\|
\nonumber\\
&\le K\eta\,\|v-\tilde{v}\|_{E_1}
\end{align*}
and $\eta>0$ is small.
Secondly, for the continuity of ${\cal S}:E_1\times E_2\to E_2$,
we let $\varsigma:=(x,y_-)$ for short and then obtain
\begin{align}
&\gamma_2^{-n}\|{\cal S}(v, w)_n(x,y_-)-{\cal S}(\tilde{v},
\tilde{w})_n(x,y_-)\|
\nonumber\\
&\le\gamma_2^{-n}\sum_{k=0}^{n-1}\|\Lambda_-^{n-k-1}\|\Big\{\|D(\pi_-
f)(v_k(\varsigma)+F^k(x))-D(\pi_- f)(\tilde{v}_k(\varsigma)+F^k(x))\|\,\|w_k(\varsigma)\|
       \nonumber\\
       &~~~~~~~+\|D(\pi_- f)(\tilde{v}_k(\varsigma)+F^k(x))\|\,\|w_k(\varsigma)-\tilde{w}_k(\varsigma)\|
                +L\|v_k(\varsigma)-\tilde{v}_k(\varsigma)\|\,\|DF^k(x)\|\Big\}
       \nonumber\\
       &~~~~+\!\gamma_2^{-n}\sum_{k=n}^{\infty}\|\Lambda_+^{n-k-1}\|\Big\{\|D(\pi_+ f)(v_k(\varsigma)\!+\!F^k(x))
             \!-\!D(\pi_+ f)(\tilde{v}_k(\varsigma)\!+\!F^k(x))\|\, \|w_k(\varsigma)\|
       \nonumber\\
       &~~~~~~~+\|D(\pi_+ f)(\tilde{v}_k(\varsigma)+F^k(x))\|\,\|w_k(\varsigma)-\tilde{w}_k(\varsigma)\|
                +L\|v_k(\varsigma)-\tilde{v}_k(\varsigma)\|\,\|DF^k(x)\|\Big\}
\nonumber\\
&\le\gamma_2^{-1}\sum_{k=0}^{n-1}\Big(\frac{\lambda_s^+}{\gamma_2}\Big)^{n-k-1}\Big\{L\gamma_2^{-k}\|v_k(\varsigma)-\tilde{v}_k(\varsigma)
       \|\big(\|w_k(\varsigma)\|+\|DF^k(x)\|\big)
       \nonumber\\
       &\hspace{4.2cm}+\eta\gamma_2^{-k}\|w_k(\varsigma)-\tilde{w}_k(\varsigma)\|\Big\}
       \nonumber\\
       &~~~~+\gamma_2^{-1}\sum_{k=n}^{\infty}\Big(\frac{\gamma_2}{\lambda_u^-}\Big)^{k+1-n}\Big\{L\gamma_2^{-k}
       \|v_k(\varsigma)-\tilde{v}_k(\varsigma)\|\big(\|w_k(\varsigma)\|+\|DF^k(x)\|\big)
       \nonumber\\
       &\hspace{4.7cm}+\eta\gamma_2^{-k}\|w_k(\varsigma)-\tilde{w}_k(\varsigma)\|\Big\}
\nonumber\\
&\le K\Big\{L\Big(\gamma_1^k\,\gamma_2^{-k}\|w_k(\varsigma)\|+\Big(\frac{\gamma_1(\lambda_u^++\eta)}{\gamma_2}\Big)^{\!k\,}\Big)
       \gamma_1^{-k}\|v_k(\varsigma)-\tilde{v}_k(\varsigma)\|
       \nonumber\\\
       &~~~~~~~~+\eta\gamma_2^{-k}\|w_k(\varsigma)-\tilde{w}_k(\varsigma)\|\Big\}
\nonumber\\
&\le
(K_1\|w\|_{E_2}+K_2)\|v-\tilde{v}\|_{E_1}+K\eta\|w-\tilde{w}\|_{E_2},
\label{xxra}
\end{align}
where $\gamma_1\in (0,1)$
because of the first inequality of (\ref{gama12}).
Moreover, setting $v=\tilde{v}$ in (\ref{xxra}), we also have
\begin{align}
\gamma_2^{-n}\|{\cal S}(v, w)_n(x,y_-)-{\cal S}(v,
\tilde{w})_n(x,y_-)\|\le K\eta\|w-\tilde{w}\|_{E_2}
\label{xxra1}
\end{align}
with small $\eta>0$. Then, we can see from (\ref{xxra}) and (\ref{xxra1}) that ${\cal S}:E_1\times E_2\to E_2$ is a continuous fiber
contraction. Therefore, (A2) is proved by Lemma \ref{lm-fbcont}.

Having (A1) and (A2), we choose a sequence
$\tilde{v}:=(\tilde{v}_n)_{n\ge 0}\in E_1$
such that $\tilde{v}_n(x,y_-)=O$ for every $n\ge 0$ and every $(x,y_-)\in \Omega$.
Clearly, $D\tilde{v}:=(D\tilde{v}_n)_{n\ge 0}\in E_2$ and each $({\cal T}\tilde{v})_n$ is
$C^1$ due to (\ref{def-T}).
Then, by the definitions (\ref{def-T})-(\ref{def-Q})
 of ${\cal T}, {\cal S}$ and ${\cal Q}$, one checks that
$
{\cal Q}(\tilde{v}, D\tilde{v})=({\cal T}\tilde{v}, D({\cal
T}\tilde{v})),
$
where $D({\cal T}\tilde{v}):=(D({\cal T}\tilde{v}_n))_{n\ge 0}$.
This enables us to prove inductively that
\begin{align}
{\cal Q}^n(\tilde{v}, D\tilde{v})=({\cal T}^n\tilde{v}, D({\cal
T}^n\tilde{v})), ~~~~~ \forall n\ge 0.
\label{Qnwv}
\end{align}
Combining (\ref{limQ}) with (\ref{Qnwv}) we get
$\lim_{n\to\infty}{\cal T}^n\tilde{v}=v_*$ and $\lim_{n\to\infty}D({\cal T}^n\tilde{v})=w_*$,
which implies that $v_*\in E_1$ such that
$Dv_*=w_*\in E_2$. Since $v_*$ is the fixed point of ${\cal T}$, it
is a solution of the Lyapunov-Perron equation (\ref{eqns-foli}).
Thus, the lemma is proved.
 \qquad$\Box$

\begin{re}
{\rm
The Lyapunov-Perron equation was also considered in \cite{ChHaTan-JDE97,CLL-JDE91,zhang-Sci00}.
In \cite{ChHaTan-JDE97,zhang-Sci00}
only the smoothness of the leaves, i.e., the smoothness of the function
$h_s(x,\cdot): X_-\to X_+$ given in (\ref{def-Ms}), is discussed.
Although the smoothness of the function $h_s: X\times X_-\to X_+$ was fully investigated
in \cite{CLL-JDE91},
their inequality (4.17) in \cite[Theorem 4.3]{CLL-JDE91} shows that
a stronger condition, i.e.,
$
\lambda_s^+\lambda_u^+<1,
$
than our (\ref{NR22}) (since $\lambda_u^->1$) was required for the existence of $C^1$ invariant foliations.
Therefore, their results are improved in
our Lemma~\ref{lm-LPeqn}.
}
\label{re-comfoli}
\end{re}


\begin{re}
{\rm
The above proof shows that the solution $(q_n)_{n\ge 0}$
of equation (\ref{eqns-foli}) is independent of $\gamma_1$ and
$\gamma_2$ because of their arbitrary choice.
In other words,
for any $\tilde{\gamma}_1$ and $\tilde{\gamma}_2$
such that
\begin{align}
\lambda_s^+<\tilde{\gamma}_1<1<\tilde{\gamma}_2<\lambda_u^-, ~~~~~
\tilde{\gamma}_1\lambda_u^+<\tilde{\gamma}_2,
\label{gamatt12}
\end{align}
the solution $(q_n)_{n\ge 0}$ obtained in Lemma \ref{lm-LPeqn} also has the inclusions:
$(q_n)_{n\ge 0}\in S_{\tilde{\gamma}_1}(\Omega,X)$ and
$(Dq_n)_{n\ge 0}\in S_{\tilde{\gamma}_2}(\Omega,{\cal L}(X\times X_-,X))$.
By the definition given in Section 2,
\begin{align}
\sup_{n\ge 0}\big\{\tilde{\gamma}_2^{-n}\|Dq_n\|_{C_b^0(\Omega,{\cal
L}(X\times X_-,X))}\big\}<\infty
\label{xyd}
\end{align}
for all $\tilde{\gamma}_2\in (\lambda_s^+\lambda_u^+, \lambda_u^-)$
due to (\ref{gamatt12}).
} \label{re-indepgama}
\end{re}


\begin{re}
{\rm The solution $q_n:\Omega \to X$ obtained in Lemma \ref{lm-LPeqn} can be extended from $\Omega$ to the whole space $X\times X_-$.
In fact, by \cite[Theorem 2.1]{ChHaTan-JDE97} we know that, for every point $(x,y_-)\in X\times X_-$,
equation (\ref{eqns-foli}) has a unique solution $(\tilde{q}_n(x,y_-))_{n\ge 0}\subset X$ such that
\begin{align*}
\sup_{n\ge 0}\big\{\gamma_1^{-n}\|\tilde{q}_n(x,y_-)\|\big\}<\infty.
\end{align*}
Note that each $\tilde{q}_n: X\times X_-\to X$ is proved to be continuous in \cite{ChHaTan-JDE97}.
By the uniqueness of $(\tilde{q}_n(x,y_-))_{n\ge 0}$ and the fact
that $(q_n)_{n\ge 0}\in S_{\gamma_1}(\Omega,X)$, we have
$
\tilde{q}_n|_{\Omega}=q_n.
$
Without loss of generality, we still let $q_n$ denote $\tilde{q}_n$ simply.
} \label{re-qnextend}
\end{re}


\section{$C^{1,\beta}$ smoothness of leaves}
\setcounter{equation}{0}

As shown in (\ref{n--0}), all $q_n$'s have the same smoothness as $q_0$ and the invariant foliation is totally determined by $q_0$.
In this section we show that $q_0:X\times X_-\to X$ obtained in Lemma \ref{lm-LPeqn} is not
only $C^1$ but also $C^{1,\beta}$ for some $\beta>0$.

\begin{lm}
Under the same conditions as in Lemma \ref{lm-LPeqn}, $q_0:X\times X_-\to X$ is $C^{1,\beta_s}$ near $O$, where
\begin{align}
\beta_s:=\frac{\log\lambda_{d}^++\log\lambda_m^+-\log\lambda_{d+1}^-}{\log\lambda_{d}^+-\log\lambda_m^+}
-\varepsilon>0
\label{def-bs}
\end{align}
for any small given $\varepsilon>0$.
\label{lm-1b}
\end{lm}

Remark that inequality (\ref{NR22}) guarantees the above-mentioned number $\beta_s$ to be positive.
In order to prove Lemma~\ref{lm-1b}, we need the following lemma to estimate the H\"{o}lder exponent of functional series.

\begin{lm}
Let $(Z,\|\cdot\|)$ be a Banach space and $\Omega\in Z$ be a small
neighborhood of the origin $\tilde{O}\in Z$. Suppose that $P_k: \Omega\to Z$ {\rm (}for all integers
$k\ge 0${\rm )} are
$C^{0,\alpha}$ mappings,  where $\alpha\in (0,1]$,
and that $\tau_1,\tau_2$ are positive constants such that
\begin{align}
\|P_k(z)\|\le M\tau_1^k, ~~~ \|P_k(z)-P_k(\tilde{z})\|\le
M\tau_2^k\|z-\tilde{z}\|^\alpha, ~~~ \forall z,\tilde{z}\in \Omega.
\label{2ineqlty}
\end{align}
If there is a number $\rho>0$ such that $\rho\tau_1<1$, then
\begin{align*}
\sum_{k=0}^{\infty}\rho^k\|P_k(z)-P_k(\tilde{z})\|\le
L\|z-\tilde{z}\|^\beta
\end{align*}
near $\tilde{O}$, where
\begin{align}
\beta:= \left\{
\begin{array}{lllll}
\alpha, &~~ \rho\tau_2<1,
\\
\alpha-\varepsilon, &~~ \rho\tau_2=1,
\\
(\log\tau_1+\log\rho)(\log\tau_1-\log\tau_2)^{-1}\,\alpha, &~~\rho\tau_2>1,
\end{array}
\right.
\label{ieq-hest}
\end{align}
defined for an arbitrarily small given $\varepsilon>0$, is a positive constant.
\label{lm-hest}
\end{lm}

{\bf Proof}. First of all, assume that $\rho\tau_2<1$. Then it is
easy to see from the second inequality of (\ref{2ineqlty}) that
\begin{align*}
\sum_{k=0}^{\infty}\rho^k\|P_k(z)-P_k(\tilde{z})\| \le
M\sum_{k=0}^{\infty}(\rho\tau_2)^k\|z-\tilde{z}\|^\alpha \le L
\|z-\tilde{z}\|^\alpha,
\end{align*}
which implies that $\beta=\alpha$.

Next, assume that $\rho\tau_2>1$. Let
$a:=\log\tau_2/\log\tau_1$ and $b:=\log\tau_2/\log\rho$
and let
\begin{align}
k_0(z,\tilde{z}):=\log\|z-\tilde{z}\|^\alpha/\log(\tau_1/\tau_2).
\label{def-k0}
\end{align}
Obviously,
\begin{align}
\tau_1^{\,a}=\tau_2, ~~~~~\rho^{\,b}=\tau_2,
~~~~~\tau_1^{k_0(z,\tilde{z})}=\tau_2^{k_0(z,\tilde{z})}\|z-\tilde{z}\|^\alpha.
\label{k0zz}
\end{align}
In what follows we let $k_0$ denote $k_0(z,\tilde{z})$ for short if
there is no confusion. Then, one computes by (\ref{k0zz}) that
\begin{align*}
\tau_1^{k_0}=\|z-\tilde{z}\|^{a(a-1)^{-1}\alpha}, ~~~
\tau_2^{k_0}=\|z-\tilde{z}\|^{(a-1)^{-1}\alpha}, ~~~
\rho^{k_0}=\|z-\tilde{z}\|^{b(a-1)^{-1}\alpha}.
\end{align*}
It follows that
\begin{align}
(\rho\tau_2)^{k_0}\|z-\tilde{z}\|^\alpha=(\rho\tau_1)^{k_0}=\|z-\tilde{z}\|^{(a+b)(a-1)^{-1}\alpha}.
\label{xbc}
\end{align}
On the other hand, $\tau_1<\tau_2$ because $\rho\tau_1<1$ and $\rho\tau_2>1$. By (\ref{def-k0}) we get
$\lim_{\|z-\tilde{z}\|\to 0}k_0(z,\tilde{z})=+\infty$, and
therefore
\begin{align}
\lim_{\|z-\tilde{z}\|\to 0}(\rho\tau_1)^{k_0}=0.
\label{zztds0}
\end{align}
Combining (\ref{xbc}) with (\ref{zztds0}), we get
$
(a+b)(a-1)^{-1}\alpha>0.
$
Putting
$$
\beta:=(a+b)(a-1)^{-1}\alpha=(\log\tau_1+\log\rho)(\log\tau_1-\log\tau_2)^{-1}\,\alpha,
$$
by (\ref{2ineqlty}) and (\ref{xbc}) we obtain
\begin{align}
\sum_{k=0}^{\infty}\rho^k\|P_k(z)-P_k(\tilde{z})\| &\le
M\sum_{k=1}^{[k_0]}(\rho\tau_2)^k\|z-\tilde{z}\|^\alpha+2M\sum_{k=[k_0]+1}^{\infty}(\rho\tau_1)^k
\nonumber\\
&\le L \|z-\tilde{z}\|^\beta.
\label{yym}
\end{align}

Finally, in the remainder case, i.e., $\rho\tau_2=1$, we note from (\ref{def-k0}) that
\begin{align*}
k_0=\log\|z-\tilde{z}\|^\alpha/\log(\tau_1/\tau_2)
=\log\|z-\tilde{z}\|^{-\alpha}/\log(\tau_2/\tau_1),
\end{align*}
where $\log(\tau_2/\tau_1)>0$ since $\tau_1<\tau_2$. Then (\ref{yym}) gives
\begin{align*}
\sum_{k=0}^{\infty}\rho^k\|P_k(z)-P_k(\tilde{z})\| &\le
M\sum_{k=1}^{[k_0]}(\rho\tau_2)^k\|z-\tilde{z}\|^\alpha+2M\sum_{k=[k_0]+1}^{\infty}(\rho\tau_1)^k
\nonumber\\
&\le Mk_0\|z-\tilde{z}\|^\alpha+K\|z-\tilde{z}\|^\alpha
\nonumber\\
&\le K(\log \|z-\tilde{z}\|^{-\alpha})
\|z-\tilde{z}\|^\alpha+K\|z-\tilde{z}\|^\alpha
\nonumber\\
&\le L\|z-\tilde{z}\|^{\alpha-\varepsilon}
\end{align*}
near $\tilde{O}$ for any given small $\varepsilon>0$.
The proof is completed. \qquad$\Box$

Now we are ready to prove Lemma \ref{lm-1b}.

{\bf Proof of Lemma \ref{lm-1b}}. Let $(q_n)_{n\ge 0}$ be the $C^1$ solution of equation (\ref{eqns-foli})
obtained in Lemma \ref{lm-LPeqn} and let
\begin{align*}
\Xi_{k}^-(x,y_-)&:=
D(\pi_-f)(q_k(x,y_-)+F^k(x))(Dq_k(x,y_-)+DF^k(x))
\nonumber\\
& ~~~~~-D(\pi_- f)(F^k(x))DF^k(x),
\end{align*}
\vspace{-0.8cm}
\begin{align*}
\Xi_{k}^+(x,y_-)
&:=
D(\pi_+f)(q_k(x,y_-)+F^k(x))(Dq_k(x,y_-)+DF^k(x))
\nonumber\\
& ~~~~~-D(\pi_+ f)(F^k(x))DF^k(x)
\end{align*}
for all integers $k\ge 0$ and all $(x,y_-),(\tilde{x},\tilde{y}_-)\in
\Omega\subset X\times X_-$. Then (\ref{eqns-foli}) gives
\begin{align}
&\gamma_2^{-n}\|Dq_n(x,y_-)-Dq_n(\tilde{x},\tilde{y}_-)\|
\nonumber\\
&\le
\gamma_2^{-1}\sum_{k=0}^{n-1}\Big(\frac{\lambda_s^+}{\gamma_2}\Big)^{n-k-1}\gamma_2^{-k}\|\Xi_{k}^-(x,y_-)-\Xi_{k}^-(\tilde{x},\tilde{y}_-)\|
\nonumber\\
&~~~~~+\gamma_2^{-1}\sum_{k=n}^{\infty}\Big(\frac{\gamma_2}{\lambda_u^-}\Big)^{k+1-n}\gamma_2^{-k}\|\Xi_{k}^+(x,y_-)-\Xi_{k}^+(\tilde{x},\tilde{y}_-)\|.
\label{xhtblt}
\end{align}
On the other hand, we claim that
\begin{align}
&\sum_{k=0}^{n-1}\Big(\frac{\lambda_s^+}{\gamma_2}\Big)^{n-k-1}\gamma_2^{-k}\|\Xi_{k}^-(x,y_-)-\Xi_{k}^-(\tilde{x},\tilde{y}_-)\|
\nonumber\\
&\le M_1
\|(x,y_-)-(\tilde{x},\tilde{y}_-)\|^{\beta_s}+M_1\eta\sup_{k\ge
0}\gamma_2^{-k}\|Dq_k(x,y_-)-Dq_k(\tilde{x},\tilde{y}_-)\|
\label{wsl-2}
\end{align}
and
\begin{align}
&\sum_{k=n}^{\infty}\Big(\frac{\gamma_2}{\lambda_u^-}\Big)^{k+1-n}\gamma_2^{-k}\|\Xi_{k}^+(x,y_-)-\Xi_{k}^+(\tilde{x},\tilde{y}_-)\|
\nonumber\\
&\le M_2
\|(x,y_-)-(\tilde{x},\tilde{y}_-)\|^{\beta_s}+M_2\eta\sup_{k\ge
0}\gamma_2^{-k}\|Dq_k(x,y_-)-Dq_k(\tilde{x},\tilde{y}_-)\|,
\label{wsl-1}
\end{align}
where $\beta_s>0$ is
given in (\ref{def-bs}). If the claimed (\ref{wsl-2}) and (\ref{wsl-1}) are both true then
(\ref{xhtblt}) gives
\begin{align*}
\sup_{n\ge
0}\gamma_2^{-n}\|Dq_n(x,y_-)\!-\!Dq_n(\tilde{x},\tilde{y}_-)\|
\le\,& M\|(x,y_-)-(\tilde{x},\tilde{y}_-)\|^{\beta_s}
\\
&+M\eta\sup_{k\ge
0}\gamma_2^{-k}\|Dq_k(x,y_-)\!-\!Dq_k(\tilde{x},\tilde{y}_-)\|,
\end{align*}
implying that
\begin{align*}
\sup_{n\ge
0}\gamma_2^{-n}\|Dq_n(x,y_-)\!-\!Dq_n(\tilde{x},\tilde{y}_-)\|\le
M(1-M\eta)^{-1}\|(x,y_-)-(\tilde{x},\tilde{y}_-)\|^{\beta_s}.
\end{align*}
Hence,
\begin{align*}
\|Dq_0(x,y_-)-Dq_0(\tilde{x},\tilde{y}_-)\|\le
L\|(x,y_-)-(\tilde{x},\tilde{y}_-)\|^{\beta_s}
\end{align*}
because $\eta$ is small, and the lemma is proved.

We only present the proof of (\ref{wsl-1}) because the proof of (\ref{wsl-2}) is similar.
By the definition of $\Xi_{k}^+(x,y_-)$, given in the beginning of the proof, and the fact that $\gamma_2/\lambda_u^-< 1$, we get
\begin{align}
&\sum_{k=n}^{\infty}\Big(\frac{\gamma_2}{\lambda_u^-}\Big)^{k+1-n}\gamma_2^{-k}\|\Xi_{k}^+(x,y_-)-\Xi_{k}^+(\tilde{x},\tilde{y}_-)\|
\nonumber\\
&\le
\sum_{k=n}^{\infty}\{\Gamma_{1k}(x,y_-,\tilde{x},\tilde{y}_-)+\Gamma_{2k}(x,y_-,\tilde{x},\tilde{y}_-)
        +\Gamma_{3k}(x,y_-,\tilde{x},\tilde{y}_-)+\Gamma_{4k}(x,y_-,\tilde{x},\tilde{y}_-)\},
\label{xi-tt}
\end{align}
where
\begin{align*}
&\Gamma_{1k}(x,y_-,\tilde{x},\tilde{y}_-)
\\
&:=\gamma_2^{-k}\|D(\pi_+ f)(q_k(x,y_-)+F^k(x))-D(\pi_+
f)(q_k(\tilde{x},\tilde{y}_-)+F^k(\tilde{x}))\|\,\|Dq_k(x,y_-)\|,
\end{align*}
\vspace{-0.8cm}
\begin{align*}
&\hspace{-0.9cm}\Gamma_{2k}(x,y_-,\tilde{x},\tilde{y}_-)
\\
&\hspace{-0.9cm}:=\Big(\frac{\gamma_2}{\lambda_u^-}\Big)^{k+1-n}\gamma_2^{-k}\|D(\pi_+
f)(q_k(\tilde{x},\tilde{y}_-)+F^k(\tilde{x}))\|
           \,\|Dq_k(x,y_-)-Dq_k(\tilde{x},\tilde{y}_-)\|,
\end{align*}
\vspace{-0.8cm}
\begin{align*}
&\hspace{-1.8cm}\Gamma_{3k}(x,y_-,\tilde{x},\tilde{y}_-)
\\
&\hspace{-1.8cm}:=\gamma_2^{-k}\|\{D(\pi_+
f)(q_k(x,y_-)+F^k(x))-D(\pi_+ f)(F^k(x))\}
\\
&\hspace{-2cm}~~~~~~~~~~~-\{D(\pi_+
f)(q_k(\tilde{x},\tilde{y}_-)+F^k(\tilde{x}))-D(\pi_+
f)(F^k(\tilde{x}))\}\|\,\|DF^k(x)\|,
\end{align*}
\vspace{-0.8cm}
\begin{align*}
&\hspace{-0.5cm}\Gamma_{4k}(x,y_-,\tilde{x},\tilde{y}_-)
\\
&\hspace{-0.5cm}:=\gamma_2^{-k}\|D(\pi_+
f)(q_k(\tilde{x},\tilde{y}_-)+F^k(\tilde{x}))-D(\pi_+
f)(F^k(\tilde{x}))\|\,\|DF^k(x)-DF^k(\tilde{x})\|.
\end{align*}
On the other hand, one can put
\begin{align*}
\gamma_1:=\lambda_s^++\delta, ~~~~~
\gamma_2:=\lambda_u^--\delta, ~~~~~
\tilde{\gamma}_2:=(\lambda_s^++\delta)\lambda_u^+
\end{align*}
for sufficiently small $\delta>0$ such that (\ref{gama12}) holds and
$\tilde{\gamma}_2\in(\lambda_s^+\lambda_u^+,\gamma_2)\subset(\lambda_s^+\lambda_u^+,\lambda_u^-)$
because of the inequality (\ref{NR22}) given in Lemma \ref{lm-LPeqn}.
Then we have the following estimates for the sums of $\Gamma_{1k}, ..., \Gamma_{4k}$ in (\ref{xi-tt}).

Concerning $\Gamma_{1k}$, we see from (\ref{xyd}) that
\begin{align*}
&\Gamma_{1k}(x,y_-,\tilde{x},\tilde{y}_-)
\\
&\le K \Big(\frac{\tilde{\gamma}_2}{\gamma_2}\Big)^k\|D(\pi_+
f)(q_k(x,y_-)+F^k(x))-D(\pi_+
f)(q_k(\tilde{x},\tilde{y}_-)+F^k(\tilde{x}))\|.
\end{align*}
Moreover, (\ref{FFbumpR2}) gives
$
\|D(\pi_+ f)(q_k(x,y_-)+F^k(x))\|\le \eta
$
and
\begin{align*}
&\|D(\pi_+ f)(q_k(x,y_-)+F^k(x))-D(\pi_+
f)(q_k(\tilde{x},\tilde{y}_-)+F^k(\tilde{x}))\|
\\
&\le
L\|q_k(x,y_-)-q_k(\tilde{x},\tilde{y}_-)\|+L\|F^k(x)-F^k(\tilde{x})\|
\\
&\le L\sup_{\varsigma\in \Omega}\|Dq_k(\varsigma)\|\,\|(x,y_-)-(\tilde{x},\tilde{y}_-)\|
      +L\sup_{\xi\in U}\|DF^k(\xi)\|\,\|x-\tilde{x}\|
\\
&\le M_1\gamma_2^k\|(x,y_-)-(\tilde{x},\tilde{y}_-)\|
     +M_2(\lambda_u^++\eta)^k\|x-\tilde{x}\|
\\
&\le M(\lambda_u^++\eta)^k\|(x,y_-)-(\tilde{x},\tilde{y}_-)\|
\end{align*}
since
$\gamma_2<\lambda_u^++\eta$.
By Lemma \ref{lm-hest}, where choosing
\begin{align*}
\alpha=1,~~~ \tau_1=1,~~~ \tau_2=\lambda_u^++\eta,~~~
\rho=\tilde{\gamma}_2/\gamma_2=(\lambda_s^++\delta)\lambda_u^+/(\lambda_u^--\delta)
\end{align*}
and noticing that $\rho\tau_1=\tilde{\gamma}_2/\gamma_2<1$,
one concludes that
\begin{align}
\sum_{k=n}^{\infty}\Gamma_{1k}(x,y_-,\tilde{x},\tilde{y}_-)\le
L\|(x,y_-)-(\tilde{x},\tilde{y}_-)\|^{\varpi_1},
\label{D1est}
\end{align}
where
\begin{align*}
\varpi_1&:=\min\Big\{1-\varepsilon,
~\frac{\log(\lambda_u^--\delta)-\log(\lambda_s^++\delta)-\log\lambda_u^+}{\log(\lambda_u^++\eta)}\Big\}
\\
&~=\min\Big\{1-\varepsilon,
~\frac{\log(\lambda_{d+1}^--2\delta)-\log(\lambda_d^++2\delta)-\log(\lambda_m^++\delta)}
        {\log(\lambda_m^++\delta+\eta)}\Big\}
\\
&~=\min\Big\{1-\varepsilon,
~\frac{\log\lambda_{d+1}^--\log\lambda_d^+-\log\lambda_m^+}{\log\lambda_m^+}
-\varepsilon\Big\}
\end{align*}
due to (\ref{ieq-hest}) for small number $\varepsilon>0$ depending on $\delta$ and $\eta$.

Concerning $\Gamma_{2k}$, we have
\begin{align}
\sum_{k=n}^{\infty}\Gamma_{2k}(x,y_-,\tilde{x},\tilde{y}_-)
&\le\sum_{k=n}^{\infty}\Big(\frac{\gamma_2}{\lambda_u^-}\Big)^{k+1-n}\eta\sup_{k\ge
0}\gamma_2^{-k}\|Dq_k(x,y_-)-Dq_k(\tilde{x},\tilde{y}_-)\|
\nonumber\\
&\le M_2\eta\sup_{k\ge
0}\gamma_2^{-k}\|Dq_k(x,y_-)-Dq_k(\tilde{x},\tilde{y}_-)\|.
\label{D2est}
\end{align}

The estimation for $\Gamma_{3k}$ is similar to $\Gamma_{1k}$.
We easily see that
\begin{align*}
\Gamma_{3k}(x,y_-,\tilde{x},\tilde{y}_-)
\le\,&\Big(\frac{\lambda_u^++\eta}{\gamma_2}\Big)^{k}\|\{D(\pi_+
f)(q_k(x,y_-)\!+\!F^k(x))\!-\!D(\pi_+ f)(F^k(x))\}
\\
&-\{D(\pi_+ f)(q_k(\tilde{x},\tilde{y}_-)+F^k(\tilde{x}))-D(\pi_+
f)(F^k(\tilde{x}))\}\|.
\end{align*}
Note that
\begin{align*}
\|D(\pi_+ f)(q_k(x,y_-)+F^k(x))-D(\pi_+ f)(F^k(x))\|\le L \|q_k(x,y_-)\|\le M \gamma_1^k
\end{align*}
and
\begin{align*}
&\|\{D(\pi_+ f)(q_k(x,y_-)+F^k(x))-D(\pi_+ f)(F^k(x))\}
\\
&~~~~~~-\{D(\pi_+f)(q_k(\tilde{x},\tilde{y}_-)+F^k(\tilde{x}))-D(\pi_+
f)(F^k(\tilde{x}))\}\|
\\
&\le M(\lambda_u^++\eta)^k\|(x,y_-)-(\tilde{x},\tilde{y}_-)\|.
\end{align*}
By Lemma \ref{lm-hest}, where choosing
\begin{align*}
\alpha=1,~~~~~ \tau_1=\gamma_1,~~~~~ \tau_2=\lambda_u^++\eta,~~~~~
\rho=(\lambda_u^++\eta)/\gamma_2
\end{align*}
and noticing that
$\rho\tau_1=\gamma_1(\lambda_u^++\eta)/\gamma_2<1 $  by (\ref{gama12}) and $\rho\tau_2=(\lambda_u^++\eta)^2/\gamma_2>1$,
one concludes that
\begin{align}
\sum_{k=n}^{\infty}\Gamma_{3k}(x,y_-,\tilde{x},\tilde{y}_-)\le
L\|(x,y_-)-(\tilde{x},\tilde{y}_-)\|^{\varpi_3},
\label{D3est}
\end{align}
where
\begin{align*}
\varpi_3&:=\frac{\log(\lambda_s^++\delta)+\log(\lambda_u^++\eta)
      -\log(\lambda_u^--\delta)}{\log(\lambda_s^++\delta)-\log(\lambda_u^++\eta)}
\\
&=\frac{\log(\lambda_d^++2\delta)+\log(\lambda_m^++\delta+\eta)-\log(\lambda_{d+1}^--2\delta)}
        {\log(\lambda_d^++2\delta)-\log(\lambda_m^++\delta+\eta)}
\\
&=\frac{\log\lambda_{d}^++\log\lambda_m^+-\log\lambda_{d+1}^-}{\log\lambda_{d}^+-\log\lambda_m^+}
-\varepsilon
\end{align*}
for small number $\varepsilon>0$ depending on $\delta,\eta>0$.
Here one can verify that $\beta$ is decreasing with respect to the small numbers $\delta$ and $\eta$.

Similarly, we have
\begin{align*}
\Gamma_{4k}(x,y_-,\tilde{x},\tilde{y}_-) &\le
L\gamma_2^{-k}\|q_k(\tilde{x},\tilde{y}_-)\|\,\|DF^k(x)-DF^k(\tilde{x})\|
\\
&\le
LK\Big(\frac{\gamma_1}{\gamma_2}\Big)^{k}\|DF^k(x)-DF^k(\tilde{x})\|.
\end{align*}
According to (\ref{FFbumpR2}), we get
$
\|DF^k(x)\|\le
(\lambda_u^++\eta)^k
$
and
\begin{align*}
&\|DF^k(x)-DF^k(\tilde{x})\|
\nonumber\\
&=
\Big\|\prod_{i=1}^{k}DF(F^{k-i}(x))-\prod_{i=1}^{k}DF(F^{k-i}(\tilde{x}))\Big\|
\nonumber\\
&\le\prod_{i=1}^{k}\|DF(F^{k-1}(\tilde{x}))\|\cdots\|DF(F^{k-i+1}(\tilde{x}))\|
                          \,\|DF(F^{k-i}(x))-DF(F^{k-i}(\tilde{x}))\|
\nonumber\\
&~~~~~~~~~~\cdot\|DF(F^{k-i-1}(x))\|\cdots\|DF(x)\|
\nonumber\\
&\le
\prod_{i=1}^{k}\bigg(\frac{(\lambda_u^++\eta)^k}{\|DF(F^{k-i}(x))\|}
         \|DF(F^{k-i}(x))-DF(F^{k-i}(\tilde{x}))\|
\bigg)
\nonumber\\
&\le
\prod_{i=1}^{k}\bigg(\frac{(\lambda_u^++\eta)^k}{\|DF(F^{k-i}(x))\|}
         L\sup_{\xi\in U}\|DF^{k-i}(\xi)\|\,\|x-\tilde{x}\|
\bigg)
\nonumber\\
&\le
\prod_{i=1}^{k}\bigg(\frac{(\lambda_u^++\eta)^k}{\|DF(F^{k-i}(x))\|}
\, L(\lambda_u^++\eta)^{k-i}\|x-\tilde{x}\|\bigg)
\nonumber\\
&\le K(\lambda_u^++\eta)^{2k}\|x-\tilde{x}\|.
\end{align*}
By Lemma \ref{lm-hest}, where choosing
\begin{align*}
\alpha=1,~~~ \tau_1=\lambda_u^++\eta,~~~
\tau_2=(\lambda_u^++\eta)^2,~~~
\rho=\gamma_1/\gamma_2=(\lambda_s^++\delta)/(\lambda_u^--\delta)
\end{align*}
and noticing that $\rho\tau_1=\gamma_1(\lambda_u^++\eta)/\gamma_2<1$,
one concludes that
\begin{align}
\sum_{k=n}^{\infty}\Gamma_{4k}(x,y_-,\tilde{x},\tilde{y}_-) &\le
L\|(x,y_-)-(\tilde{x},\tilde{y}_-)\|^{\varpi_4}
\label{D4est}
\end{align}
where
\begin{align*}
\varpi_4&:=\min\Big\{1-\varepsilon,
~\frac{\log(\lambda_u^--\delta)-\log(\lambda_s^++\delta)-\log(\lambda_u^++\eta)}{\log(\lambda_u^++\eta)}\Big\}
\\
&~=\min\Big\{1-\varepsilon,
~\frac{\log\lambda_{d+1}^--\log\lambda_d^+-\log\lambda_m^+}{\log\lambda_m^+}
-\varepsilon\Big\}.
\end{align*}
Therefore, (\ref{wsl-1}) is proved by (\ref{xi-tt})-(\ref{D4est}) because one checks that
\begin{align*}
\beta_s&=\min\{\varpi_1,1,\varpi_3,\varpi_4\}
\\
&=\min\Big\{1-\varepsilon,
~\frac{\log\lambda_{d+1}^--\log\lambda_d^+-\log\lambda_m^+}{\log\lambda_m^+}
-\varepsilon,
\frac{\log\lambda_{d}^++\log\lambda_m^+-\log\lambda_{d+1}^-}{\log\lambda_{d}^+-\log\lambda_m^+}
-\varepsilon\Big\}
\\
&=\frac{\log\lambda_{d}^++\log\lambda_m^+-\log\lambda_{d+1}^-}{\log\lambda_{d}^+-\log\lambda_m^+}
-\varepsilon.
\end{align*}
The proof is completed. \qquad$\Box$


\section{Partial linearization of contractions}
\setcounter{equation}{0}

As indicated in the Introduction, $F$ will be decomposed for a reduction to linearization of a contraction and
linearization of an expansion.
Even only for a contraction $F$, the procedure of linearization will be divided into several steps
according to the decomposition of the spectrum and in each step
we linearize $F$ along a corresponding band in the decomposition.
In this section we discuss the linearization along a band of the spectrum,
called a {\it partial linearization}.

Assume that (\ref{spec-decomp})-(\ref{spec-order}) hold. For the sake of convenience,
we use another version of the Spectral Decomposition Theorem given in \cite{El-JFA01} to
assume that $X=X_1\times\cdots \times X_m$ with $\Lambda$-invariant subspaces $X_i$'s, that is,
\begin{align*}
\Lambda={\rm diag} (\Lambda_1,...,\Lambda_m),
\end{align*}
where $\Lambda_i\in {\cal L}(X_i, X_i)$ such that
$\sigma(\Lambda_i)=\sigma_i$ for all $i=1,...,m$.
Our strategy is to realize the linearization of $F$ with a sequence of
diffeomorphisms $\Phi_{i}$, each of which makes a conjugacy between
$F_{i}$ and $F_{{i}-1}$, i.e., the partial linearization. Here $F_{i}$'s have the
following forms: $F_0=\Lambda$,
\begin{align}
F_{i}(x)= \left(\begin{array}{lllllll}
\Lambda_1x_1+f_{1{i}}(x_1,...,x_m)
\\
~~~~~~~~~~~~\vdots
\\
\Lambda_{i} x_{i}+f_{{i}{i}}(x_1,...,x_m)
\\
\Lambda_{{i}+1}x_{{i}+1}
\\
~~~~~~~~~~~~\vdots
\\
\Lambda_{m}x_{m}
\end{array}\right),
~~~~~ \forall {i}=1,...,m-1,
\label{F1-k}
\end{align}
and $F_m=F$, where $x=(x_1,...,x_m)\in X_1\times\cdots \times X_m$
and $f_{j{i}}:X\to X_j$ for all $j=1,...,{i}$. The transferring
process is shown in the diagram:
\begin{align}
F=F_m \xrightarrow{~\Phi_{m}~} \cdots \xrightarrow{\Phi_{{i}+1}}
F_{i} \xrightarrow{~\Phi_{{i}}~} F_{{i}-1}
\xrightarrow{\Phi_{{i}-1}} \cdots \xrightarrow{~\Phi_{1}~}
F_0=\Lambda.
\label{transff}
\end{align}

In order to investigate these partial linearizations of $F$, for a
given integer $\ell\in \{2,...,m-1\}$ we use the notations
\begin{align*}
{\cal U}:=X_1\times\cdots \times X_{\ell-1}, ~~~~ {\cal V}:=X_\ell,
~~~~ {\cal W}:=X_{\ell+1}\times\cdots \times X_m
\end{align*}
and let
\begin{align*}
{A}:={\rm diag}(\Lambda_1,...,\Lambda_{\ell-1}), ~~~~
{B}:=\Lambda_\ell, ~~~~ {C}:={\rm
diag}(\Lambda_{\ell+1},...,\Lambda_{m})
\end{align*}
accordingly. Then $F_\ell$ defined in (\ref{F1-k}) with
$i=\ell$ can be rewritten as
\begin{align}
F_\ell(x)=(Au+ f_a(u,v,w), Bv+f_b(u,v,w), Cw),
\label{F-induc}
\end{align}
where $x:=(u,v,w)\in{\cal U}\times{\cal V}\times{\cal W}$, $f_a:X\to
{\cal U}$ and $f_b:X\to {\cal V}$. By \cite[Theorem 5]{R-S-JDE04},
for any sufficiently small $\delta>0$, we can appropriately choose
equivalent norms in ${\cal U},{\cal V}$ and ${\cal W}$ separately
such that
\begin{align}
\|A\|\le \|B\|\le \mu_\ell^+, ~~~ \|B^{-1}\|\le 1/\mu_\ell^-, ~~~
\|C\|\le \mu_m^+, ~~~ \|C^{-1}\|\le 1/\mu_{\ell+1}^-,
\label{def-delta}
\end{align}
where $\mu_i^+:=\lambda_i^+ +\delta$, $\mu_i^-:=\lambda_i^-
-\delta$, $i=\ell,\ell+1,m$, and use the norm
$\|x\|:=\|u\|+\|v\|+\|w\|$ for $x\in X$. Obviously,
\begin{align}
0<\mu_\ell^-<\mu_\ell^+<\mu_{\ell+1}^-<\mu_m^+<1.
\label{mumumu}
\end{align}
We further assume that $f_a$ and $f_b$ given in (\ref{F-induc}) are
both $C^1$ such that
\begin{align}
f_a(O)=f_b(O)=0, ~~~~~  Df_a(O)=Df_b(O)=0
\label{fafb}
\end{align}
and the inequalities
\begin{align}
& \|\partial_{(u,v)}f_a(x)-\partial_{(u,v)}f_a(\tilde{x})\|\le M_1\|x-\tilde{x}\|,
\label{hold-l1}
\\
& \|\partial_{(u,v)}f_b(x)-\partial_{(u,v)}f_b(\tilde{x})\|\le M_1\|x-\tilde{x}\|,
\label{hold-l2}
\\
& \|\partial_{w}f_a(u,v,w)-\partial_{w}f_a(\tilde{u},\tilde{v},w)\|
        \le M_1\|(u,v)-(\tilde{u},\tilde{v})\|,
\label{hold-l3}
\\
&\|\partial_{w}f_a(u,v,w)-\partial_{w}f_a(u,v,\tilde{w})\|
        \le M_1\|w-\tilde{w}\|^{\beta_{\ell+1}},
\label{hold-l4}
\\
& \|\partial_{w}f_b(u,v,w)-\partial_{w}f_b(\tilde{u},\tilde{v},w)\|
        \le M_1\|(u,v)-(\tilde{u},\tilde{v})\|,
\label{hold-l5}
\\
& \|\partial_{w}f_b(u,v,w)-\partial_{w}f_b(u,v,\tilde{w})\|
        \le M_1\|w-\tilde{w}\|^{\beta_{\ell+1}}
\label{hold-l6}
\end{align}
hold in $U$, where $\beta_{\ell+1}\in(0,1]$. Inequalities
(\ref{hold-l1})-(\ref{hold-l6}) obviously imply that $F_\ell$ is
$C^{1,\beta_{\ell+1}}$ in $U$. Our lemma on partial linearization is the
following:

\begin{lm}
Suppose that the
$C^1$ mapping ${F_\ell}:X\to X$, given in {\rm (\ref{F-induc})},
satisfies {\rm (\ref{fafb})-(\ref{hold-l6})} and that inequality {\rm
(\ref{NR-conB})} with $i=\ell$ holds, i.e.,
\begin{align}
\lambda_\ell^+/\lambda_\ell^-<(\lambda_m^+)^{-1}.
\label{NR-conB-lll}
\end{align}
Then there exists a $C^1$ diffeomorphism
$\Phi_\ell: U\to X$ such that
\begin{align}
F_{\ell-1}(x):=\Phi_\ell\circ
F_\ell\circ\Phi_\ell^{-1}(x)=(Au+{g}(u,v,w), Bv, Cw),
\label{Fl-1-induc}
\end{align}
where $g:U\to {\cal U}$ is $C^1$ such that $ g(O)=0$ and $Dg(O)=0. $
Moreover, $\Phi_\ell$ is $C^{1,\beta_{\ell}}$ and the inequalities
\begin{align}
& \|\partial_{(u,v)}g(x)-\partial_{(u,v)}g(\tilde{x})\|\le M_2\|x-\tilde{x}\|,
\label{hold-l+11}
\\
& \|\partial_{w}g(u,v,w)-\partial_{w}g(\tilde{u},\tilde{v},w)\|\le
M_2 \|(u-\tilde{u},v-\tilde{v})\|,
\label{hold-l+12}
\\
&\|\partial_{w}g(u,v,w)-\partial_{w}g(u,v,\tilde{w})\|\le M_2
\|w-\tilde{w}\|^{\beta_{\ell}}
\label{hold-l+13}
\end{align}
hold in $U$, where $\beta_{\ell}\in (0,\beta_{\ell+1}]$ is defined by
\begin{align}
\left\{
\begin{array}{ll}
\beta_{\ell}:=\min\Big\{\zeta-\varepsilon,
~\frac{\log\lambda_\ell^++\log\lambda_m^+-\log\lambda_\ell^-}
     {\log\lambda_\ell^+-\zeta\log\lambda_m^+}\,\zeta-\varepsilon\Big\},
\\
\vspace{-0.1cm}
\\
\zeta:=\min\{\beta_{\ell+1},~\log\lambda_\ell^+/\log\lambda_{\ell+1}^-
-1-\varepsilon\}
\end{array}
\right.
\label{def-betal}
\end{align}
for any small given $\varepsilon>0$.
\label{lm-induc}
\end{lm}

Before proving Theorem \ref{lm-induc}, we need a lemma on invariant manifolds of $F_\ell$.

\begin{lm}
Suppose that
${F_\ell}:X\to X$, given in {\rm (\ref{F-induc})},  is $C^{1,\zeta}$ such that
both {\rm (\ref{def-delta})}
and {\rm (\ref{fafb})} hold with the inequality
\begin{align}
0<\zeta<\log\mu_\ell^+/\log\mu_{\ell+1}^- -1.
\label{<zata<}
\end{align}
Then $F_\ell$ has a local $C^{1,\zeta}$ invariant manifold
\begin{align*}
\Gamma:=\{(u,v,w)\in U: u=h_1(w),v=h_2(w)\},
\end{align*}
where both $h_1: U\cap {\cal W} \to {\cal U}$ and $h_2: U\cap {\cal W} \to {\cal V}$
are $C^{1,\zeta}$ such that $h_1(0)=h_2(0)=0$ and $Dh_1(0)=Dh_2(0)=0$.
\label{lm-inv}
\end{lm}

Remark that the number $\zeta$ given in (\ref{<zata<}) is well defined since $\log\mu_\ell^+/\log\mu_{\ell+1}^->1$
by (\ref{mumumu}).

{\bf Proof}.
This lemma is actually a corollary of \cite[Theorem 2.1]{Llave-Wayne-MZ1995} and
the fourth remark given in \cite{Llave-Wayne-MZ1995}.
In fact,
by (\ref{def-delta}) and (\ref{<zata<}), there exists a sufficiently
small constant $\epsilon>0$ such that
\begin{align*}
\|(F_\ell^{-1})'(O)|_{{\cal W}}\|=\|C^{-1}\|\le 1/\mu_{\ell+1}^-<a
~~ {\rm and} ~~
\|F_\ell'(O)|_{{\cal U}\times{\cal V}}\|=\|B\|\le \mu_\ell^+<a^{-(1+\zeta)},
\end{align*}
where $a:=(\mu_{\ell+1}^- -\epsilon)^{-1}$.
\quad $\Box$

{\bf Proof of Lemma \ref{lm-induc}}. Firstly, we simplify
$F_\ell$ near $O$ by straightening up its an invariant manifold. For this purpose,
we note that $F_\ell$ is
$C^{1,\beta_{\ell+1}}$ in $U$, as indicated in
(\ref{hold-l1})-(\ref{hold-l6}), such that (\ref{def-delta}) holds.
Moreover, for any small given $\varepsilon>0$, one can always find small $\delta>0$ such that
$
\log\lambda_\ell^+/\log\lambda_{\ell+1}^-
-1-\varepsilon<\log\mu_\ell^+/\log\mu_{\ell+1}^--1
$
by the definitions of those $\mu_i$'s given just below (\ref{def-delta}).
Thus, putting
\begin{align}
\zeta:=\min\{\beta_{\ell+1},~\log\lambda_\ell^+/\log\lambda_{\ell+1}^-
-1-\varepsilon\}
\label{zatatwo}
\end{align}
and according to Lemma \ref{lm-inv}, we conclude that
$F_\ell$ has a local $C^{1,\zeta}$ invariant manifold
\begin{align*}
\Gamma=\{(u,v,w)\in U: u=h_1(w),v=h_2(w)\},
\end{align*}
where both $h_1: U\cap {\cal W} \to {\cal U}$ and $h_2: U\cap {\cal
W} \to {\cal V}$ are $C^{1,\zeta}$ such that $h_1(0)=h_2(0)=0$ and
$Dh_1(0)=Dh_2(0)=0$.
This enables us to define a diffeomorphism $\Theta:U\to X$ by
\begin{align}
\Theta(x):=(u-h_1(w), v-h_2(w), w)
\label{def-theta}
\end{align}
and define a new mapping ${\widetilde{F}_\ell}:U\to X$ by
\begin{align}
\hspace{-0.2cm}
{\widetilde{F}_\ell}(x):=\Theta\circ F_\ell\circ\Theta^{-1}(x)
=\left(\begin{array}{ll}
\pi_uF_\ell(u+h_1(w),v+h_2(w),w)-h_1({{C}}w)
\\
\pi_vF_\ell(u+h_1(w),v+h_2(w),w)-h_2({{C}}w)
\\
{{C}}w
\end{array}\right),
\label{Fww}
\end{align}
where $\pi_u$ and $\pi_v$ denote the projections onto subspaces
${\cal U}$ and ${\cal V}$ respectively. Then it is easy to verify
that
\begin{align}
\pi_u{\widetilde{F}_\ell}(0,0,w)=\pi_v{\widetilde{F}_\ell}(0,0,w)=0,
~~~~~ \forall w\in U\cap {\cal W}.
\label{F00w=0}
\end{align}
Moreover, according to (\ref{fafb})-(\ref{hold-l6}) and (\ref{Fww}),
one checks by straightforward computation that
$\widetilde{F}_\ell(O)=O$, $D\widetilde{F}_\ell(O)=\Lambda$ and
\begin{align}
&
\|\partial_{(u,v)}(\pi_u{\widetilde{F}_\ell})(x)-\partial_{(u,v)}(\pi_u{\widetilde{F}_\ell})(\tilde{x})\|
     \le M_1\|x-\tilde{x}\|,
\label{hold-l1ww}
\\
&
\|\partial_{(u,v)}(\pi_v{\widetilde{F}_\ell})(x)-\partial_{(u,v)}(\pi_v{\widetilde{F}_\ell})(\tilde{x})\|
     \le M_1\|x-\tilde{x}\|,
\label{hold-l2ww}
\\
&
\|\partial_{w}(\pi_u{\widetilde{F}_\ell})(u,v,w)-\partial_{w}(\pi_u{\widetilde{F}_\ell})(\tilde{u},\tilde{v},w)\|
     \le M_1 \|(u,v)-(\tilde{u},\tilde{v})\|,
\label{hold-l3ww}
\\
&\|\partial_{w}(\pi_u{\widetilde{F}_\ell})(u,v,w)-\partial_{w}(\pi_u{\widetilde{F}_\ell})(u,v,\tilde{w})\|
     \le M_1\|w-\tilde{w}\|^{\zeta},
\label{hold-l4ww}
\\
&
\|\partial_{w}(\pi_v{\widetilde{F}_\ell})(u,v,w)-\partial_{w}(\pi_v{\widetilde{F}_\ell})(\tilde{u},\tilde{v},w)\|
     \le M_1\|(u,v)-(\tilde{u},\tilde{v})\|,
\label{hold-l5ww}
\\
&\|\partial_{w}(\pi_v{\widetilde{F}_\ell})(u,v,w)-\partial_{w}(\pi_v{\widetilde{F}_\ell})(u,v,\tilde{w})\|
     \le M_1\|w-\tilde{w}\|^{\zeta}.
\label{hold-l6ww}
\end{align}

In order to prove the first part of the lemma, i.e., the part before
``Moreover", we show the uniform convergence of the sequence
$(B^{-k}D(\pi_v{\widetilde{F}_\ell}^k))_{k\in\mathbb{N}}$ in $U$. If
the uniform convergence is true,
then the sequence
$(B^{-k}\pi_v{\widetilde{F}_\ell}^k)_{k\in\mathbb{N}}$ is also
convergent uniformly in $U$ by the Mean Value Theorem in
Banach spaces. Thus, the limit
\begin{align*}
\widetilde{\Psi}_\ell:=\lim_{k\to\infty}B^{-k}\pi_v{\widetilde{F}_\ell}^k
\end{align*}
gives a $C^1$ mapping from $U$ into ${\cal V}$ and we can define
\begin{align}
\Psi_\ell:=\widetilde{\Psi}_\ell\circ\Theta,
\label{tetFlta}
\end{align}
which satisfies that
\begin{align}
\Psi_\ell\circ F_\ell
&=\widetilde{\Psi}_\ell\circ\Theta\circ
F_\ell
   =\widetilde{\Psi}_\ell\circ \widetilde{F}_\ell\circ\Theta
   =\lim_{k\to\infty}B^{-k}\pi_v{\widetilde{F}_\ell}^{k+1}\circ\Theta
\nonumber\\
&=B\lim_{k\to\infty}B^{-(k+1)}\pi_v{\widetilde{F}_\ell}^{k+1}\circ\Theta=B\widetilde{\Psi}_\ell\circ\Theta
   =B\Psi_\ell.
\label{psil}
\end{align}
Note that $D\Psi_\ell(O)=(0,{\rm id}_{{\cal
V}},0)$ by (\ref{def-theta}) and (\ref{tetFlta}), where ${\rm
id}_{{\cal V}}$ denotes the identity mapping in ${\cal V}$. Then we
can further define a $C^1$ diffeomorphism $\Phi_\ell:U\to X$ by
\begin{align}
\Phi_\ell(x):=(u,\Psi_\ell(u,v,w),w)
\label{philll}
\end{align}
and, by the well-known Inverse Function Theorem, its inverse is of the form
\begin{align}
\Phi_\ell^{-1}(x):=(u,\widehat{\Psi}_\ell(u,v,w),w),
\label{philll-}
\end{align}
where $\widehat{\Psi}_\ell:U\to{\cal V}$ has the same regularity as
$\Psi_\ell$. It follows from (\ref{psil})-(\ref{philll-}) that
\begin{align}
\Phi_\ell\circ F_\ell\circ\Phi_\ell^{-1}(x)
&=(\pi_u
F_{\ell}\circ\Phi_\ell^{-1}(x), \Psi_\ell\circ
F_\ell\circ\Phi_\ell^{-1}(x), Cw)
\nonumber\\
&=(Au+{g}(x), B\Psi_\ell\circ\Phi_\ell^{-1}(x), Cw)
\nonumber\\
&=(Au+{g}(x), Bv, Cw),
\label{hhh}
\end{align}
where
\begin{align}
{g}(x):=\pi_u
F_{\ell}\circ\Phi_\ell^{-1}(x)-Au=f_a(u,\widehat{\Psi}_\ell(u,v,w),w)
\label{Fl-1Fl}
\end{align}
by (\ref{F-induc}).
In fact,
$\Psi_\ell\circ\Phi_\ell^{-1}(x)=v$ in the third equality of
(\ref{hhh}) because
\begin{align*}
\Phi_\ell\circ\Phi_\ell^{-1}(x)=(u,\Psi_\ell\circ\Phi_\ell^{-1}(x),w)=(u,v,w).
\end{align*}
One checks that $g:U\to {\cal U}$ defined in (\ref{Fl-1Fl})
 is $C^1$ such that $ g(O)=0$ and $Dg(O)=0
$ since $\Phi_\ell(O)=O$ and $D\Phi_\ell(O)={\rm id}$, the identity
mapping in $X$. Therefore $\Phi_\ell$ is the desired diffeomorphism
in the lemma.

In what follows, we prove the uniform convergence of the above-mentioned sequence
$(B^{-k}D(\pi_v{\widetilde{F}_\ell}^k))_{k\in\mathbb{N}}$. For
convenience, put
\begin{eqnarray}
D{\widetilde{F}_\ell}^k :=\left(\begin{array}{cccc} \hat{a}_{k} &
\hat{b}_{k} & \hat{c}_{k}
\\
a_{k} & b_{k} & c_{k}
\\
0 & 0 & C^k
\end{array}\right),
\label{abcn} ~~~~~ \forall k\in \mathbb{N},
\label{def-Fn'}
\end{eqnarray}
where $ \hat{a}_k:=\partial_{u}(\pi_u \widetilde{F}_\ell^k),
\hat{b}_k:=\partial_{v}(\pi_u \widetilde{F}_\ell^k),
\hat{c}_k:=\partial_{w}(\pi_u \widetilde{F}_\ell^k),
a_k:=\partial_{u}(\pi_v \widetilde{F}_\ell^k),
b_k:=\partial_{v}(\pi_v \widetilde{F}_\ell^k) $ and
$c_k:=\partial_{w}(\pi_v \widetilde{F}_\ell^k)$. Obviously,
\begin{equation}
\begin{array}{l}
a_{k+1}(x)=a_1({\widetilde{F}_\ell}^k(x))\hat{a}_k(x)+b_1({\widetilde{F}_\ell}^k(x))a_k(x),
\\
b_{k+1}(x)=a_1({\widetilde{F}_\ell}^k(x))\hat{b}_k(x)+b_1({\widetilde{F}_\ell}^k(x))b_k(x),
\\
c_{k+1}(x)=a_1({\widetilde{F}_\ell}^k(x))\hat{c}_k(x)+b_1({\widetilde{F}_\ell}^k(x))c_k(x)+c_1({\widetilde{F}_\ell}^k(x))C^k
\end{array}
\label{dtgs}
\end{equation}
because $D(\widetilde{F}_\ell^{k+1})(x)=DF_\ell(\widetilde{F}_\ell^k(x))D\widetilde{F}_\ell^k(x)$.
The following lemma is necessary.
\begin{lm}
The inequalities
\begin{align}
& \max\{\|\hat{a}_k(x)\|,~ \|a_k(x)\|,~
\|\hat{b}_k(x)\|,~ \|b_k(x)\|, ~ \|\hat{c}_k(x)\|,~ \|c_k(x)\|\}
      \le K(\mu_\ell^+)^{k},
\label{estim1}
\\
& \max\{\|a_1({\widetilde{F}_\ell}^k(x))\|,~
\|b_1({\widetilde{F}_\ell}^k(x))-{{B}}\|\}
     \le K(\mu_m^+)^{k}, ~ \|c_1({\widetilde{F}_\ell}^k(x))\|
     \le K(\mu_\ell^+)^{k}
\label{estim2}
\end{align}
hold in $U$ for all $k\in\mathbb{N}$.
\label{lm-123}
\end{lm}

We postpone the proof of this lemma to next section but just apply
it to continue our proof of the theorem. According to (\ref{dtgs}),
\begin{align}
&\|B^{-(k+1)}a_{k+1}(x)-B^{-k}a_{k}(x)\|\le {\cal S}_{1k}(x),
\label{esti-xi1}
\\
& \|B^{-(k+1)}b_{k+1}(x)-B^{-k}b_{k}(x)\|\le {\cal S}_{2k}(x),
\label{esti-xi2}
\\
& \|B^{-(k+1)}c_{k+1}(x)-B^{-k}c_{k}(x)\|\le {\cal S}_{3k}(x),
\label{esti-xi3}
\end{align}
where
\begin{align*}
&{\cal S}_{1k}(x)
:=\|B^{-(k+1)}\|\Big(\|a_1({\widetilde{F}_\ell}^k(x))\|\,\|\hat{a}_k(x)\|+\|b_1({\widetilde{F}_\ell}^k(x))-B\|\,\|a_k(x)\|\Big),
\\
&{\cal S}_{2k}(x)
:=\|B^{-(k+1)}\|\Big(\|a_1({\widetilde{F}_\ell}^k(x))\|\,\|\hat{b}_k(x)\|+\|b_1({\widetilde{F}_\ell}^k(x))-B\|\,\|b_k(x)\|\Big),
\\
&{\cal S}_{3k}(x)
:=\|B^{-(k+1)}\|\Big(\|a_1({\widetilde{F}_\ell^k(x)})\|\,\|\hat{c}_k(x)\|+\|b_1({\widetilde{F}_\ell^k(x)})-B\|\,\|c_k(x)\|
\\
&~~~~~~~~~~~~~~~~~~~~~~~~~~~~~
   +\|c_1({\widetilde{F}_\ell^k(x)})\|\,\|C^k\|\Big).
\end{align*}
From (\ref{def-delta}) and inequalities (\ref{estim1})-
(\ref{estim2}), given in Lemma~\ref{lm-123}, we see that
\begin{align}
\max\{{\cal S}_{1k}(x),{\cal S}_{2k}(x),{\cal S}_{3k}(x)\}\le
K\eta^k,
\label{S123}
\end{align}
where
\begin{align}
\eta:=\mu_m^+\mu_\ell^+/\mu_\ell^-.
\label{eta12}
\end{align}
By (\ref{NR-conB-lll}) one checks that $0<\eta<1$ since $\mu_i$'s can be
chosen arbitrarily close to $\lambda_i$'s, as indicated below
(\ref{def-delta}). Therefore, the sequence
$(B^{-k}D(\pi_v{\widetilde{F}_\ell}^k))_{k\in\mathbb{N}}$ is
uniformly convergent in $U$ because
\begin{align*}
&\lim_{k\to\infty}B^{-k}D(\pi_v{\widetilde{F}_\ell}^k)
\nonumber\\
&=\sum_{k=1}^{\infty}\{B^{-(k+1)}D(\pi_v \widetilde{F}_\ell^{k+1})
-B^{-k}D(\pi_v
\widetilde{F}_\ell^{k})\}+B^{-1}D(\pi_v{\widetilde{F}_\ell})
\end{align*}
and
\begin{align}
&\|B^{-(k+1)}D(\pi_v \widetilde{F}_\ell^{k+1})(x)-B^{-k}D(\pi_v
\widetilde{F}_\ell^{k})(x)\|
\nonumber\\
&\le \max\Big\{\|B^{-(k+1)}a_{k+1}(x)-B^{-k}a_{k}(x)\|,\|B^{-(k+1)}b_{k+1}(x)-B^{-k}b_{k}(x)\|,
\nonumber\\
&\hspace{1.7cm}\|B^{-(k+1)}c_{k+1}(x)-B^{-k}c_{k}(x)\|\Big\}
\nonumber\\
&\le\max\{{\cal S}_{1k}(x),{\cal S}_{2k}(x),{\cal S}_{3k}(x)\} \le K\eta^k
\label{BvF}
\end{align}
by (\ref{esti-xi1})-(\ref{S123}). Thus the proof of the first part
of the lemma is completed.

In what follows, we prove the second part of the lemma, i.e., the
part after ``Moreover". In order to prove
(\ref{hold-l+11})-(\ref{hold-l+13}), by (\ref{Fl-1Fl}) we see that
it suffices to investigate the regularity of $\Psi_\ell$
defined in (\ref{tetFlta}) because $\widehat{\Psi}_\ell$ has the
same regularity as $\Psi_\ell$, as mentioned just below
(\ref{philll-}). For this purpose, we need to notice that
$\Psi_\ell=\widetilde{\Psi}_\ell\circ\Theta$ by (\ref{tetFlta}),
where $
\widetilde{\Psi}_\ell:=\lim_{k\to\infty}B^{-k}\pi_v{\widetilde{F}_\ell}^k.
$ We next claim that
\begin{align}
&
\|\partial_{(u,v)}\widetilde{\Psi}_\ell(x)-\partial_{(u,v)}\widetilde{\Psi}_\ell(\tilde{x})\|
\nonumber\\
&
=\|\lim_{k\to\infty}B^{-k}\{\partial_{(u,v)}(\pi_v{\widetilde{F}_\ell}^k)(x)-\partial_{(u,v)}(\pi_v{\widetilde{F}_\ell}^k)(\tilde{x})\}\|
     \le M_2\|x-\tilde{x}\|,
~~~~~~~~~~~~\label{estim-xy1}
\end{align}
\begin{align}
&
\|\partial_{w}\widetilde{\Psi}_\ell(u,v,w)-\partial_{w}\widetilde{\Psi}_\ell(\tilde{u},\tilde{v},w)\|
\nonumber\\
&
=\|\lim_{k\to\infty}B^{-k}\{\partial_{w}(\pi_v{\widetilde{F}_\ell}^k)(u,v,w)
    -\partial_{w}(\pi_v{\widetilde{F}_\ell}^k)(\tilde{u},\tilde{v},w)\}\|
    \le M_2\|(u,v)-(\tilde{u},\tilde{v})\|
~~~ \label{estim-xy2}
\end{align}
and
\begin{align}
&
\|\partial_{w}\widetilde{\Psi}_\ell(u,v,w)-\partial_{w}\widetilde{\Psi}_\ell(u,v,\tilde{w})\|
\nonumber\\
&
=\|\lim_{k\to\infty}B^{-k}\{\partial_w(\pi_v{\widetilde{F}_\ell}^k)(u,v,w)-\partial_w(\pi_v{\widetilde{F}_\ell}^k)(u,v,\tilde{w})\}\|
     \le M_2\|w-\tilde{w}\|^{\beta_{\ell}}
 \label{estim-xy4}
\end{align}
for a number $\beta_{\ell}\in (0,\zeta]$ defined by
\begin{align}
\beta_{\ell}:=\min\Big\{\zeta-\varepsilon,
~\frac{\log\lambda_\ell^++\log\lambda_m^+-\log\lambda_\ell^-}
     {\log\lambda_\ell^+-\zeta\log\lambda_m^+}\,\zeta-\varepsilon\Big\},
\label{betal-1ym}
\end{align}
where $\zeta$ is
given in (\ref{zatatwo}). If the claimed
(\ref{estim-xy1})-(\ref{estim-xy4}) are all true, then by
(\ref{def-theta}) and (\ref{tetFlta}) one can check directly that
$\Psi_\ell$ satisfies the same inequalities as in
(\ref{estim-xy1})-(\ref{estim-xy4}) and so does the mapping
$\widehat{\Psi}_\ell$. Thus we can prove that $\Phi_\ell$ is
$C^{1,\beta_{\ell}}$ by (\ref{philll}) and that inequalities
(\ref{hold-l+11})-(\ref{hold-l+13}) hold by (\ref{hold-l1}),
(\ref{hold-l3}) and (\ref{hold-l4}) because
\begin{align*}
\partial_u{g}(x)&=\partial_uf_a(u,\widehat{\Psi}_\ell(x),w)+\partial_vf_a(u,\widehat{\Psi}_\ell(x),w)\cdot\partial_u\widehat{\Psi}_\ell(x),
\\
\partial_v{g}(x)&=\partial_vf_a(u,\widehat{\Psi}_\ell(x),w)\cdot\partial_v\widehat{\Psi}_\ell(x),
\\
\partial_w{g}(x)&=\partial_vf_a(u,\widehat{\Psi}_\ell(x),w)\cdot\partial_w\widehat{\Psi}_\ell(x)+\partial_wf_a(u,\widehat{\Psi}_\ell(x),w),
\end{align*}
as known from (\ref{Fl-1Fl}). This completes the proof of the second
part of the theorem.

Now we complementarily prove the claimed
(\ref{estim-xy1})-(\ref{estim-xy4}). We will use the following
lemma.

\begin{lm}
The inequalities
\begin{align}
& \max\{\|\hat{a}_{k}(x)-\hat{a}_{k}(\tilde{x})\|, \|a_{k}(x)-a_{k}(\tilde{x})\|,
    \|\hat{b}_{k}(x)-\hat{b}_{k}(\tilde{x})\|, \|b_{k}(x)-b_{k}(\tilde{x})\|\}
\nonumber\\
&~~~~~~\le L(\mu_\ell^+)^k\|x-\tilde{x}\|,
\label{estim1-D}
\\
\nonumber\\
&
\max\{\|\hat{c}_{k}(u,v,w)\!-\!\hat{c}_{k}(\tilde{u},\tilde{v},w)\|,
     \|c_{k}(u,v,w)-c_{k}(\tilde{u},\tilde{v},w)\|\}
\nonumber\\
&~~~~~~ \le L(\mu_\ell^+)^k\|(u,v)-(\tilde{u},\tilde{v})\|,
\label{estim11-D}
\\
\nonumber\\
& \max\{\|\hat{c}_{k}(u,v,w)-\hat{c}_{k}(u,v,\tilde{w})\|,
      \|c_{k}(u,v,w)-c_{k}(u,v,\tilde{w})\|\}
\nonumber\\
&~~~~~~ \le L(\mu_m^+)^k\|w-\tilde{w}\|^{\zeta},
\label{estim14-D}
\\
\nonumber\\
&
\max\{\|a_1({\widetilde{F}_\ell}^k(x))-a_1({\widetilde{F}_\ell}^k(\tilde{x}))\|,
    \|b_1({\widetilde{F}_\ell}^k(x))-b_1({\widetilde{F}_\ell}^k(\tilde{x}))\|\}
\nonumber\\
&~~~~~~\le L(\mu_m^+)^{k}\|x-\tilde{x}\|,
\label{estim4-D}
\\
\nonumber\\
&
\|c_1({\widetilde{F}_\ell}^k(u,v,w))-c_1({\widetilde{F}_\ell}^k(\tilde{u},\tilde{v},w))\|
     \le L(\mu_\ell^+)^{k}\|(u,v)-(\tilde{u},\tilde{v})\|,
\label{estim13-D}
\\
\nonumber\\
&
\|c_1({\widetilde{F}_\ell}^k(u,v,w))-c_1({\widetilde{F}_\ell}^k(u,v,\tilde{w}))\|
      \le L(\mu_m^+)^{k\zeta}\|w-\tilde{w}\|^{\zeta}
\label{estim16-D}
\end{align}
hold in $U$ for all $k\in\mathbb{N}$.
\label{lm-dc}
\end{lm}

We still leave the proof of Lemma~\ref{lm-dc} to next section. Using
(\ref{dtgs}) again,
\begin{align*}
&\|B^{-(k+1)}(a_{k+1}(x)-a_{k+1}(\tilde{x}))-B^{-k}(a_k(x)-a_k(\tilde{x}))\|\le{\cal
D}_{1k}(x,\tilde{x}),
\\
&\|B^{-(k+1)}(b_{k+1}(x)-b_{k+1}(\tilde{x}))-B^{-k}(b_k(x)-b_k(\tilde{x}))\|\le{\cal
D}_{2k}(x,\tilde{x}),
\\
&\|B^{-(k+1)}(c_{k+1}(x)-c_{k+1}(\tilde{x}))-B^{-k}(c_k(x)-c_k(\tilde{x}))\|\le{\cal
D}_{3k}(x,\tilde{x}),
\end{align*}
where
\begin{align*}
& {\cal D}_{1k}(x,\tilde{x})
\\
&
:=\|B^{-(k+1)}\|\Big(\|a_1({\widetilde{F}_\ell}^k(x))-a_1({\widetilde{F}_\ell}^k(\tilde{x}))\|\,\|\hat{a}_k(x)\|
        +\|a_1({\widetilde{F}_\ell}^k(\tilde{x}))\|\,\|\hat{a}_k(x)-\hat{a}_k(\tilde{x})\|
\\
&~~~~~+\|b_1({\widetilde{F}_\ell}^k(x))-b_1({\widetilde{F}_\ell}^k(\tilde{x}))\|\,\|a_k(x)\|
        +\|b_1({\widetilde{F}_\ell}^k(\tilde{x}))-B\|\,\|a_k(x)-a_k(\tilde{x})\|\Big),
\\
\\
& {\cal D}_{2k}(x,\tilde{x})
\\
&
:=\|B^{-(k+1)}\|\Big(\|a_1({\widetilde{F}_\ell}^k(x))-a_1({\widetilde{F}_\ell}^k(\tilde{x}))\|\,\|\hat{b}_k(x)\|
        +\|a_1({\widetilde{F}_\ell}^k(\tilde{x}))\|\,\|\hat{b}_k(x)-\hat{b}_k(\tilde{x})\|
~~~~~~~~~~~
\\
&~~~~~+\!\|b_1({\widetilde{F}_\ell}^k(x))-b_1({\widetilde{F}_\ell}^k(\tilde{x}))\|\,\|b_k(x)\|
        +\|b_1({\widetilde{F}_\ell}^k(\tilde{x}))-B\|\,\|b_k(x)-b_k(\tilde{x})\|\Big),
\\
\\
& {\cal D}_{3k}(x,\tilde{x})
\\
&
:=\|B^{-(k+1)}\|\Big(\|a_1({\widetilde{F}_\ell}^k(x))-a_1({\widetilde{F}_\ell}^k(\tilde{x}))\|\,\|\hat{c}_k(x)\|
        +\|a_1({\widetilde{F}_\ell}^k(\tilde{x}))\|\,\|\hat{c}_k(x)-\hat{c}_k(\tilde{x})\|
~~~~~~~~~~~
\\
&~~~~~+\|b_1({\widetilde{F}_\ell}^k(x))-b_1({\widetilde{F}_\ell}^k(\tilde{x}))\|\,\|c_k(x)\|
        +\|b_1({\widetilde{F}_\ell}^k(\tilde{x}))-B\|\,\|c_k(x)-c_k(\tilde{x})\|
\\
&~~~~~+\|c_1({\widetilde{F}_\ell}^k(x))-c_1({\widetilde{F}_\ell}^k(\tilde{x}))\|\,\|C^k\|\Big).
\end{align*}
It follows from Lemmas \ref{lm-123} and \ref{lm-dc}
that
\begin{align*}
\max\{{\cal D}_{1k}(x,\tilde{x}),~
{\cal D}_{2k}(x,\tilde{x})\}\le
\frac{4KL}{\mu_\ell^-}\eta^k\|x-\tilde{x}\|,
\end{align*}
where $\eta\in (0,1)$ is given in (\ref{eta12}). Then, an analogous argument as (\ref{BvF}) gives
\begin{align*}
&
\Big\|\lim_{k\to\infty}B^{-k}\Big(\partial_{(u,v)}(\pi_v{\widetilde{F}_\ell}^k)(x)-\partial_{(u,v)}(\pi_v{\widetilde{F}_\ell}^k)(\tilde{x})\Big)\Big\|
\\
& \le\sum_{k=1}^{\infty}\max\{{\cal D}_{1k}(x,\tilde{x}),{\cal D}_{2k}(x,\tilde{x})\}
    +\|B^{-1}\|\,\|\partial_{(u,v)}(\pi_v{\widetilde{F}_\ell})(x)-\partial_{(u,v)}(\pi_v{\widetilde{F}_\ell})(\tilde{x})\|
\\
& \le M_2\|x-\tilde{x}\|
\end{align*}
by (\ref{hold-l2ww}). This proves (\ref{estim-xy1}). Similarly,
Lemmas \ref{lm-123} and \ref{lm-dc}
yield
\begin{align*}
{\cal D}_{3k}((u,v,w),(\tilde{u},\tilde{v},w))
\le \frac{5KL}{\mu_\ell^-}\eta^k
\|(u,v)-(\tilde{u},\tilde{v})\|,
\end{align*}
which implies that
\begin{align*}
&
\Big\|\lim_{k\to\infty}B^{-k}\Big(\partial_{w}(\pi_v{\widetilde{F}_\ell}^k)(u,v,w)
-\partial_{w}(\pi_v{\widetilde{F}_\ell}^k)(\tilde{u},\tilde{v},w)\Big)\Big\|
\\
& \le\sum_{k=1}^{\infty}{\cal
D}_{3k}((u,v,w),(\tilde{u},\tilde{v},w))
      +\|B^{-1}\|\,\|\partial_{w}(\pi_v{\widetilde{F}_\ell})(u,v,w)-\partial_{w}(\pi_v{\widetilde{F}_\ell})(\tilde{u},\tilde{v},w)\|
\\
& \le M_2\|(u,v)-(\tilde{u},\tilde{v})\|
\end{align*}
by (\ref{hold-l5ww}). This proves (\ref{estim-xy2}).

In order to prove the claimed (\ref{estim-xy4}), we need the following
lemma.
\begin{lm}
The inequalities
\begin{align}
&
\sum_{k=1}^{\infty}\|B^{-(k+1)}\|\,\|a_1({\widetilde{F}_\ell}^k(x))\|\,\|\hat{c}_k(u,v,w)-\hat{c}_k(u,v,\tilde{w})\|
        \le M\|w-\tilde{w}\|^{\beta_{\ell}},
\label{se1}
\\
&
\sum_{k=1}^{\infty}\|B^{-(k+1)}\|\,\|b_1({\widetilde{F}_\ell}^k(x))\!-\!B\|\,\|c_k(u,v,w)-c_k(u,v,\tilde{w})\|
        \le M\|w-\tilde{w}\|^{\beta_{\ell}},
\label{se2}
\\
&
\sum_{k=1}^{\infty}\|B^{-(k+1)}\|\,\|c_1({\widetilde{F}_\ell}^k(u,v,w))-c_1({\widetilde{F}_\ell}^k(u,v,\tilde{w}))\|\,\|C^k\|
        \le M\|w-\tilde{w}\|^{\beta_{\ell}}
\label{se3}
\end{align}
hold in $U$ for all $k\in\mathbb{N}$, where $\beta_{\ell}\in
(0,\zeta]$ is the number defined by {\rm (\ref{betal-1ym})}.
\label{lm-w}
\end{lm}

The proof of Lemma~\ref{lm-w} will be given in next section. By
(\ref{hold-l6ww}) and Lemmas \ref{lm-123}-\ref{lm-w},
we get
\begin{align*}
&
\|\lim_{k\to\infty}B^{-k}(\partial_{w}(\pi_v{\widetilde{F}_\ell}^k)(u,v,w)-\partial_{w}(\pi_v{\widetilde{F}_\ell}^k)(u,v,\tilde{w}))\|
\\
&
\le\sum_{k=1}^{\infty}{\cal
D}_{3k}((u,v,w),(u,v,\tilde{w}))
   +\|B^{-1}\|\,\|\partial_{w}(\pi_v{\widetilde{F}_\ell})(u,v,w)-\partial_{w}(\pi_v{\widetilde{F}_\ell})(u,v,\tilde{w})\|
\\
&
\le\frac{2KL}{\mu_\ell^-}\sum_{k=1}^{\infty}\eta^k\|w-\tilde{w}\|+3M\|w-\tilde{w}\|^{\beta_{\ell}}
      +\frac{M_1}{\mu_\ell^-}\|w-\tilde{w}\|^{\zeta}
      \le M_2\|w-\tilde{w}\|^{\beta_{\ell}}
\end{align*}
because $\beta_{\ell}\le \zeta\le 1$, which proves
(\ref{estim-xy4}). The proof is completed. \quad$\Box$

\begin{re}
{\rm The inequality $\|c_1({\widetilde{F}_\ell}^k(x))\|\le
K(\mu_\ell^+)^{k}$, given in (\ref{estim2}),
has a subtle
difference from the corresponding estimates of $a_1({\widetilde{F}_\ell}^k(x))$ and
$b_1({\widetilde{F}_\ell}^k(x))-B$ but
plays a critical role in the proof of Lemma \ref{lm-induc}.
Indeed, the number $\mu_\ell^+$ in (\ref{eta12})
comes from the inequality $\|c_1({\widetilde{F}_\ell}^k(x))\|\le
K(\mu_\ell^+)^{k}$, as indicated in the definition of ${\cal S}_{3k}$ given above (\ref{S123}).
The number $\eta$
in (\ref{eta12}) cannot be smaller than 1 if we replace $\mu_\ell^+$ with a larger number $\mu_m^+$.
}
\label{re-ineq}
\end{re}

\begin{re}
{\rm In this section the main result is Lemma \ref{lm-induc}, which
completes the the partial linearization in the inductive
transferring process (\ref{transff}). $\Phi_\ell$ is actually a
$C^{1,\beta_\ell}$ solution of a functional equation deduced from
the conjugacy $\Phi_\ell\circ F_\ell=F_{\ell-1}\circ \Phi_\ell$ and
its derivative $D\Phi_\ell$ needs to satisfy some additional
properties (\ref{hold-l+11})-(\ref{hold-l+13}) so that the induction
can be applied. In the proof of Lemma~\ref{lm-induc} we did not
employ the Contraction Principle, which demands a complicated
functional space of $C^{1,\beta_\ell}$ mappings
and many complicated estimates in a special norm, but
constructed a sequence of mappings to approximate $\Phi_\ell$
in the $C^1$ norm simply. After the $C^1$ solution $\Phi_\ell$ is
obtained, we verified that $\Phi_\ell$ is $C^{1,\beta_\ell}$ and
satisfies (\ref{hold-l+11})-(\ref{hold-l+13}). This method provides
a chance to find better estimates for those exponents $\beta_\ell$
and therefore we can obtain better estimates for
$\beta$. } \label{re-meth}
\end{re}


\section{Proofs of Lemmas}
\setcounter{equation}{0}

In this section we prove Lemmas \ref{lm-123}-\ref{lm-w}.

{\bf Proof of Lemma \ref{lm-123}}. From
(\ref{hold-l1ww})-(\ref{hold-l6ww}) we see that $\widetilde{F}_\ell$
is $C^{1,\zeta}$, where $\zeta>0$ is given in (\ref{zatatwo}). Then
\begin{align}
\|D\widetilde{F}_\ell^k(x)\|
&\le\prod_{i=0}^{k-1}\|D\widetilde{F}_\ell(\widetilde{F}_\ell^i(x))\|
                 \le\prod_{i=0}^{k-1}(\|\Lambda\|+M_1\|\widetilde{F}_\ell^i(x)\|^\zeta)
\nonumber\\
&\le\prod_{i=0}^{k-1}(\mu_m^++M_1\sup_{\xi\in
U}\|D\widetilde{F}_\ell^i(\xi)\|^\zeta)
      \le\prod_{i=0}^{k-1}(\mu_m^++M_1\prod_{j=0}^{i-1}\sup_{\xi\in
U}\|D\widetilde{F}_\ell(\widetilde{F}_\ell^j(\xi))\|^\zeta)
\nonumber\\
&\le\prod_{i=0}^{k-1}(\mu_m^++M_1(\mu_m^++\eta)^{i
\zeta})
                 \le K_1(\mu_m^+)^k
\label{fn'<uk}
\end{align}
for all $k\in \mathbb{N}$ and all $x\in U$, where $\eta>0$ is
so small that
$
\mu_m^++\eta<1
$
and
$
K_1\ge\prod_{i=0}^{\infty}(1+\frac{M_1}{\mu_m^+}(\mu_m^++\eta)^{i
\zeta})\in (0,\infty).
$
Let
\begin{align}
Q_k(x):= \left(\begin{array}{cccc} \hat{a}_{k}(x) & \hat{b}_{k}(x)
\\
a_{k}(x) & b_{k}(x)
\end{array}\right),
~~~~~ \forall k\in \mathbb{N},
\label{def-qn}
\end{align}
where $\hat{a}_{k},\hat{b}_{k},a_{k}$ and $b_{k}$ are given in
(\ref{abcn}). Then,
\begin{align}
Q_{k+1}(x)=Q_1(\widetilde{F}_\ell^k(x))Q_k(x)
\label{Qn+1Qn}
\end{align}
because $D(\widetilde{F}_\ell^{k+1})(x)=DF_\ell(\widetilde{F}_\ell^k(x))D\widetilde{F}_\ell^k(x)$.
Note that $Q_1$ is $C^{1,1}$ by (\ref{hold-l1ww}) and
(\ref{hold-l2ww}). Moreover, $Q_1(O)={\rm diag}(A,B)$ since $D\widetilde{F}_\ell(O)={\rm diag}(A,B,C)$. Then
(\ref{fn'<uk}) yields
\begin{align}
\|Q_1(\widetilde{F}_\ell^k(x))-{\rm diag}(A,B)\|
&\le
M_1\|\widetilde{F}_\ell^k(x)\| \le M_1\sup_{\xi\in
U}\|D\widetilde{F}_\ell^k(\xi)\|\|x\|
\nonumber\\
&\le M_1K_1(\mu_m^+)^{k}\|x\|,
\label{Q111}
\end{align}
which proves the first inequality given in (\ref{estim2}).

Combining (\ref{Qn+1Qn}) with (\ref{Q111}), we get
\begin{align*}
\|Q_{i+1}(x)\|
&\le\|Q_1(\widetilde{F}_\ell^i(x))\|\,\|Q_{i}(x)\|
      \le(\|{\rm diag}(A,B)\|+M_1K_1(\mu_m^+)^i)\|Q_{i}(x)\|
\\
&\le(\mu_\ell^++M_1K_1(\mu_m^+)^i)\|Q_{i}(x)\|, ~~~~~
\forall i\in\mathbb{N},
\end{align*}
and therefore
\begin{align}
\|Q_{k}(x)\|\le \|Q_1(x)\|\prod_{i=1}^{k-1}(\mu_\ell^+
+M_1K_1(\mu_m^+)^i)\le K_2(\mu_\ell^+)^{k}, ~~~~~\forall
k\in\mathbb{N},
\label{qnnn}
\end{align}
where $ K_2\ge \sup_{\xi\in
U}\|Q_1(\xi)\|(\mu_\ell^+)^{-1}\prod_{i=1}^{\infty}(1+\frac{M_1K_1}{\mu_\ell^+}(\mu_m^+)^i)\in(0,\infty).
$
According to (\ref{F00w=0}),
\begin{align*}
\pi_u{\widetilde{F}_\ell^k}(0,0,w)=\pi_v{\widetilde{F}_\ell^k}(0,0,w)=0,
~~~
\partial_{w}(\pi_u{\widetilde{F}_\ell})(0,0,w)=\partial_{w}(\pi_v{\widetilde{F}_\ell})(0,0,w)=0
\end{align*}
for all $k\in\mathbb{N}$ and $w\in U\cap {\cal W}$. Thus,
(\ref{qnnn}) gives
\begin{align}
&\|(\pi_u{\widetilde{F}_\ell^k}(x),\pi_v{\widetilde{F}_\ell^k}(x))\|
\nonumber\\
&=\|(\pi_u{\widetilde{F}_\ell^k}(u,v,w),\pi_v{\widetilde{F}_\ell^k}(u,v,w))
     -(\pi_u{\widetilde{F}_\ell^k}(0,0,w),\pi_v{\widetilde{F}_\ell^k}(0,0,w))\|
\nonumber\\
&\le \sup_{\xi\in U}\|Q_k(\xi)\|\,\|(u,v)\|\le
K_2(\mu_\ell^+)^{k}\|(u,v)\|
\label{Fuv}
\end{align}
and (\ref{hold-l5ww}) gives
\begin{align}
\|c_1(x)\|&=\|\partial_{w}(\pi_v{\widetilde{F}_\ell})(x)\|
\nonumber\\
&=\|\partial_{w}(\pi_v{\widetilde{F}_\ell})(u,v,w)-\partial_{w}(\pi_v{\widetilde{F}_\ell})(0,0,w)\|
\le M_1\|(u,v)\|. \label{c1uv}
\end{align}
It follows from (\ref{Fuv}) and (\ref{c1uv}) that
\begin{align}
\|c_1({\widetilde{F}_\ell}^k(x))\| \le
M_1\|(\pi_u{\widetilde{F}_\ell}(x),\pi_v{\widetilde{F}_\ell}(x))\|
\le M_1K_2(\mu_\ell^+)^{k}\|(u,v)\|
\label{c1fn1}
\end{align}
for all $k\in\mathbb{N}$. This proves the second inequality given in
(\ref{estim2}).

Since (\ref{def-qn}) and (\ref{qnnn}) yield
\begin{align}
\max\{\|\hat{a}_k(x)\|,~ \|a_k(x)\|,~ \|\hat{b}_k(x)\|, ~
\|b_k(x)\|\}\le K_2(\mu_\ell^+)^{k},
\label{aabb}
\end{align}
it remains to prove that
\begin{align}
\max\{\|\hat{c}_k(x)\|,~ \|c_k(x)\|\}\le K(\mu_\ell^+)^{k}.
\label{cc}
\end{align}
Similarly to (\ref{c1fn1}),
we get
\begin{align}
\|\hat{c}_1({\widetilde{F}_\ell}^k(x))\| \le
M_1\|(\pi_u{\widetilde{F}_\ell}(x),\pi_v{\widetilde{F}_\ell}(x))\|
\le M_1K_2(\mu_\ell^+)^{k}\|(u,v)\|.
\label{c1fn2}
\end{align}
By (\ref{def-Fn'}) and (\ref{def-qn}) we have
\begin{align}
\!\!\!\!(\hat{c}_{i+1}(x),c_{i+1}(x))^T =
Q_1({\widetilde{F}_\ell}^i(x))(\hat{c}_{i}(x),c_{i}(x))^T
   +(\hat{c}_1({\widetilde{F}_\ell}^i(x))C^i,c_1({\widetilde{F}_\ell}^i(x))C^i)^T
\label{cn+1cn}
\end{align}
for all $i\in\mathbb{N}$, where the superscript $T$ denotes the
transpose of vectors. Then
\begin{align*}
&\|(\hat{c}_{i+1}(x),c_{i+1}(x))^T\|
\\
&
\le\|Q_1({\widetilde{F}_\ell}^i(x))\|\,\|(\hat{c}_{i}(x),c_{i}(x))^T\|
     +\|(\hat{c}_1({\widetilde{F}_\ell}^i(x)),c_1({\widetilde{F}_\ell}^i(x)))^T\|\,\|C^i\|
\\
&\le (\mu_\ell^+ +M_1K_1(\mu_m^+)^{i})\|(\hat{c}_{i}(x),c_{i}(x))^T\|
     +2M_1K_2(\mu_\ell^+\mu_m^+)^i
\end{align*}
because of (\ref{Q111}), (\ref{c1fn1}) and (\ref{c1fn2}). Since
$\hat{c}_{1}(x)$ and $c_{1}(x)$ are both bounded in $U$, the
positive constants $M_1$ and $K_2$ can be chosen so large that $\|(\hat{c}_{1}(x),c_{1}(x))\|\le 2M_1K_2$. Therefore
\begin{align*}
& \|(\hat{c}_{k}(x),c_{k}(x))^T\|\le2M_1K_2\sum_{j=0}^{k-1}
     \Big\{(\mu_\ell^+\mu_m^+)^j\prod_{i=j+1}^{k-1}(\mu_\ell^+ +M_1K_1(\mu_m^+)^i)\Big\}
\end{align*}
for all $k\in\mathbb{N}$ by induction, where $\prod_{i=j+1}^{k-1}(\mu_\ell^+ +M_1K_1(\mu_m^+)^i)=1$ when $j=k-1$.
Hence we get
\begin{align}
\frac{\|(\hat{c}_{k}(x),c_{k}(x))^T\|}{(\mu_\ell^+)^{k}}
&\le
\frac{2M_1K_2}{\mu_\ell^+}\sum_{j=0}^{k-1}\Big\{(\mu_m^+)^j
      \prod_{i=j+1}^{k-1}\Big(1+\frac{M_1K_1}{\mu_\ell^+}(\mu_m^+)^i\Big)\Big\}
\nonumber\\
&\le \frac{2M_1K_2}{\mu_\ell^+}\sum_{j=0}^{\infty}(\mu_m^+)^j
      \prod_{i=1}^{\infty}\Big(1+\frac{M_1K_1}{\mu_\ell^+}(\mu_m^+)^i\Big)<\infty
\label{cncn}
\end{align}
because $\mu_m^+<1$, which proves (\ref{cc}).
Thus inequality (\ref{estim1}) is obtained from (\ref{aabb}) and (\ref{cc}). The proof is completed.
\qquad$\Box$

{\bf Proof of Lemma \ref{lm-dc}}. Recall that
$\widetilde{F}_\ell$ is $C^{1,\zeta}$ in $U$. It follows from
(\ref{fn'<uk}) that
\begin{align}
&\|D\widetilde{F}_\ell(\widetilde{F}_\ell^k(x))-D\widetilde{F}_\ell(\widetilde{F}_\ell^k(\tilde{x}))\|
\nonumber\\
&\le M_1\|\widetilde{F}_\ell^k(x)-\widetilde{F}_\ell^k(\tilde{x})\|^\zeta
\le M_1\sup_{\xi\in
U}\|D\widetilde{F}_\ell^k(\xi)\|^\zeta\|x-\tilde{x}\|^\zeta
\nonumber\\
&\le M_1K_1^\zeta(\mu_m^+)^{k\zeta}\|x-\tilde{x}\|^\zeta
\label{xy1}
\end{align}
for all $k\in \mathbb{N}$ and all $x,\tilde{x}\in U$, which implies
(\ref{estim16-D}). Similarly, since $Q_1$ is $C^{1,1}$ as mentioned
below (\ref{Qn+1Qn}), we have
\begin{align}
\|Q_1(\widetilde{F}_\ell^k(x))\!-\!Q_1(\widetilde{F}_\ell^k(\tilde{x}))\|
\le M_1\|\widetilde{F}_\ell^k(x)-\widetilde{F}_\ell^k(\tilde{x})\| \le
M_1K_1(\mu_m^+)^{k}\|x-\tilde{x}\|,
\label{q1xy1}
\end{align}
which proves (\ref{estim4-D}). Combining (\ref{fn'<uk}) with
(\ref{xy1}), we get
\begin{align}
&\|D\widetilde{F}_\ell^k(x)-D\widetilde{F}_\ell^k(\tilde{x})\|
\nonumber\\
&\le
\sum_{i=1}^{k}\Bigg{(}\frac{\|D\widetilde{F}_\ell(\widetilde{F}_\ell^{k-1}(\tilde{x}))\|
                               \cdots\|D\widetilde{F}_\ell(\widetilde{F}_\ell^{k-i+1}(\tilde{x}))\|
        \,\|D\widetilde{F}_\ell(\widetilde{F}_\ell^{k-i}(x))\|\cdots\|D\widetilde{F}_\ell(x)\|}
                               {\|D\widetilde{F}_\ell(\widetilde{F}_\ell^{k-i}(x))\|}
\nonumber\\
&
~~~~~~~~~~~~\cdot\|D\widetilde{F}_\ell(\widetilde{F}_\ell^{k-i}(x))-D\widetilde{F}_\ell(\widetilde{F}_\ell^{k-i}(\tilde{x}))\|
\Bigg{)}
\nonumber\\
&\le \sum_{i=1}^{k}\Bigg{(}
\frac{K_1(\mu_m^+)^k}{\|D\widetilde{F}_\ell(\widetilde{F}_\ell^{k-i}(x))\|}
       M_1K_1^\zeta(\mu_m^+)^{(k-i)\zeta}\|x-\tilde{x}\|^{\zeta}\Bigg{)}
\nonumber\\
&\le L(\mu_m^+)^{k}\|x-\tilde{x}\|^{\zeta}.
\label{xany1}
\end{align}
Here $\|D\widetilde{F}_\ell(\widetilde{F}_\ell^{k-i}(x))\|^{-1}$ is bounded near $O$ because
$\|D\widetilde{F}_\ell(O)\|=\|\Lambda\|\ne 0$. This proves
(\ref{estim14-D}).

From (\ref{Qn+1Qn})-(\ref{qnnn}) and (\ref{q1xy1}) one sees that
\begin{align*}
& \|Q_{i+1}(x)-Q_{i+1}(\tilde{x})\|
\\
&=\|Q_1(\widetilde{F}_\ell^i(x))Q_{i}(x)-Q_1(\widetilde{F}_\ell^i(\tilde{x}))Q_{i}(\tilde{x})\|
\\
&\le\|Q_1(\widetilde{F}_\ell^i(x))-Q_1(\widetilde{F}_\ell^i(\tilde{x}))\|\,\|Q_{i}(x)\|
      +\|Q_1(\widetilde{F}_\ell^i(\tilde{x}))\|\,\|Q_{i}(x)-Q_{i}(\tilde{x})\|
\\
&\le (\mu_\ell^++M_1K_1(\mu_m^+)^i)\|Q_{i}(x)-Q_{i}(\tilde{x})\|
         +M_1K_1 K_2(\mu_\ell^+\mu_m^+)^{i}\|x-\tilde{x}\|
\end{align*}
for all $i\in\mathbb{N}$. Then the same arguments as before
(\ref{cncn}) give
\begin{align*}
\frac{\|Q_{k}(x)-Q_{k}(\tilde{x})\|}{(\mu_\ell^+)^{k}\|x-\tilde{x}\|} \le
\frac{M_1K_1
K_2}{\mu_\ell^+}\prod_{i=1}^{\infty}\Big(1+\frac{M_1K_1}{\mu_\ell^+}(\mu_m^+)^i\Big)
      \sum_{j=0}^{\infty}(\mu_m^+)^j<\infty
\end{align*}
for all $k\in\mathbb{N}$, which proves (\ref{estim1-D}).

Finally, we prove the remaining two inequalities (\ref{estim11-D})
and (\ref{estim13-D}). Noticing that $\pi_w\widetilde{F}_\ell(x)=Cw$ as seen in (\ref{def-Fn'}),
we prove (\ref{estim13-D}) because
\begin{align}
&
\|c_1({\widetilde{F}_\ell}^k(u,v,w))-c_1({\widetilde{F}_\ell}^k(\tilde{u},\tilde{v},w))\|
\nonumber\\
&
\le M_1\|(\pi_u{\widetilde{F}_\ell}^k(u,v,w),\pi_v{\widetilde{F}_\ell}^k(u,v,w))
     -(\pi_u{\widetilde{F}_\ell}^k(\tilde{u},\tilde{v},w),\pi_v{\widetilde{F}_\ell}^k(\tilde{u},\tilde{v},w))\|
\nonumber\\
&
\le M_1\sup_{\xi\in U}\|Q_k(\xi)\|\,\|(u,v)-(\tilde{u},\tilde{v})\|
     \le M_1K_2(\mu_\ell^+)^{k}\|(u,v)-(\tilde{u},\tilde{v})\|
\label{c1-c1}
\end{align}
due to (\ref{hold-l5ww}) and (\ref{qnnn}). Similarly,
\begin{align}
\|\hat{c}_1({\widetilde{F}_\ell}^k(u,v,w))-\hat{c}_1({\widetilde{F}_\ell}^k(\tilde{u},\tilde{v},w))\|
\le M_1K_2(\mu_\ell^+)^k\|(u,v)-(\tilde{u},\tilde{v})\|.
\label{c1-c1hat}
\end{align}
On the other hand, (\ref{cn+1cn}) yields
\begin{align*}
&
\|(\hat{c}_{i+1}(u,v,w),c_{i+1}(u,v,w))^T-(\hat{c}_{i+1}(\tilde{u},\tilde{v},w),c_{i+1}(\tilde{u},\tilde{v},w))^T\|
\\
& \le
\|Q_1({\widetilde{F}_\ell}^i(u,v,w))-Q_1({\widetilde{F}_\ell}^i(\tilde{u},\tilde{v},w))\|\,\|(\hat{c}_i(u,v,w),c_i(u,v,w))^T\|
\\
&~~~~+\|Q_1({\widetilde{F}_\ell}^i(\tilde{u},\tilde{v},w))\|\,\|(\hat{c}_i(u,v,w),c_i(u,v,w))^T
            -(\hat{c}_i(\tilde{u},\tilde{v},w),c_i(\tilde{u},\tilde{v},w))^T\|
\\
&~~~~+\|(\hat{c}_1({\widetilde{F}_\ell}^i(u,v,w)),c_1({\widetilde{F}_\ell}^i(u,v,w)))^T
            -(\hat{c}_1({\widetilde{F}_\ell}^i(\tilde{u},\tilde{v},w)),c_1({\widetilde{F}_\ell}^i(\tilde{u},\tilde{v},w)))^T\|
            \,\|C^i\|.
\end{align*}
Consequently, by (\ref{Q111}), (\ref{cc}) and
(\ref{q1xy1})-(\ref{c1-c1hat}),
\begin{align*}
&
\|(\hat{c}_{i+1}(u,v,w),c_{i+1}(u,v,w))^T-(\hat{c}_{i+1}(\tilde{u},\tilde{v},w),c_{i+1}(\tilde{u},\tilde{v},w))^T\|
\\
& \le
2KM_1K_1(\mu_\ell^+\mu_m^+)^{i}\|(u,v)-(\tilde{u},\tilde{v})\|
\\
     &~~~~+(\mu_\ell^+ +M_1K_1(\mu_m^+)^i)\|(\hat{c}_i(u,v,w),c_i(u,v,w))^T
           -(\hat{c}_i(\tilde{u},\tilde{v},w),c_i(\tilde{u},\tilde{v},w))^T\|
\\
     &~~~~+2M_1K_2(\mu_\ell^+\mu_m^+)^{i}\|(u,v)-(\tilde{u},\tilde{v})\|
\\
& = (\mu_\ell^+
+M_1K_1(\mu_m^+)^i)\|(\hat{c}_i(u,v,w),c_i(u,v,w))^T-(\hat{c}_i(\tilde{u},\tilde{v},w),c_i(\tilde{u},\tilde{v},w))^T\|
\\
     &~~~~+2(KM_1K_1+M_1K_2)\,(\mu_\ell^+\mu_m^+)^{i}\,\|(u,v)-(\tilde{u},\tilde{v})\|
\end{align*}
for all $i\in\mathbb{N}$. The same arguments as before (\ref{cncn})
give
\begin{align*}
&
\frac{\|((\hat{c}_k(u,v,w),c_k(u,v,w))^T-(\hat{c}_k(\tilde{u},\tilde{v},w),c_k(\tilde{u},\tilde{v},w))^T\|}
       {(\mu_\ell^+)^{k}\|(u,v)-(\tilde{u},\tilde{v})\|}
\\
& \le \frac{2(KM_1K_1+M_1K_2)}{\mu_\ell^+}
      \sum_{j=0}^{\infty}(\mu_m^+)^{j}
      \prod_{i=1}^{\infty}\Big(1+\frac{M_1K_1}{\mu_\ell^+}(\mu_m^+)^{i}\Big)<\infty
\end{align*}
for all $k\in\mathbb{N}$. This proves
(\ref{estim11-D}) and the proof is completed. \quad$\Box$

{\bf Proof of Lemma \ref{lm-w}}. For (\ref{se1}), note that $\|c_k(u,v,w)\|\le K(\mu_\ell^+)^{k}$
by (\ref{estim1}) and that
\begin{align*}
\|c_k(u,v,w)-c_k(u,v,\tilde{w})\|\le
L(\mu_m^+)^{k}\|w-\tilde{w}\|^{\zeta}
\end{align*}
by (\ref{estim14-D}). In view of Lemma \ref{lm-hest}, where choosing
\begin{align*}
\alpha=\zeta, ~~~~~\tau_1=\mu_\ell^+, ~~~~~ \tau_2=\mu_m^+, ~~~~~
\rho=\mu_m^+/\mu_\ell^-,
\end{align*}
and noticing that $\rho\tau_1=\eta<1$, we obtain
\begin{align}
&\sum_{k=1}^{\infty}\|B^{-(k+1)}\|\,\|a_1({\widetilde{F}_\ell}^k(x))\|\,\|\hat{c}_k(u,v,w)-\hat{c}_k(u,v,\tilde{w})\|
\nonumber\\
&\le\frac{K}{\mu_\ell^-}\sum_{k=1}^{\infty}
       \Big(\frac{\mu_m^+}{\mu_\ell^-}\Big)^k\|\hat{c}_k(u,v,w)-\hat{c}_k(u,v,\tilde{w})\|
\nonumber\\
&
\le M\|w-\tilde{w}\|^{p_1},
\label{se-proof1}
\end{align}
where $p_1:=\min\{\zeta-\varepsilon,\,\varpi\}$ and
\begin{align*}
\varpi:=\frac{\log\mu_\ell^++\log\mu_m^+-\log\mu_\ell^-}{\log\mu_\ell^+-\log\mu_m^+}\,\zeta
=\frac{\log(\lambda_\ell^++\delta)+\log(\lambda_m^++\delta)-\log(\lambda_\ell^--\delta)}
              {\log(\lambda_\ell^++\delta)-\log(\lambda_m^++\delta)}\,\zeta
\end{align*}
In the case of $\lambda_\ell^-<(\lambda_m^+)^2$, we understand that $\varpi$ is decreasing with respect
to the small number $\delta>0$ and therefore
\begin{align*}
\varpi=\frac{\log\lambda_\ell^++\log\lambda_m^+-\log\lambda_\ell^-}
     {\log\lambda_\ell^+-\log\lambda_m^+}\,\zeta-\varepsilon>0
\end{align*}
for small number $\varepsilon>0$ depending on $\delta$. In the case of $\lambda_\ell^-\ge(\lambda_m^+)^2$
we have $\varpi>\zeta-\varepsilon$. Hence, in either case
\begin{align*}
p_1=\min\{\zeta-\varepsilon,\,\varpi\}
=\min\Big\{\zeta-\varepsilon,\,\frac{\log\lambda_\ell^++\log\lambda_m^+-\log\lambda_\ell^-}
     {\log\lambda_\ell^+-\log\lambda_m^+}\,\zeta-\varepsilon\Big\}.
\end{align*}
For (\ref{se2}), using the same arguments we get
\begin{align}
\sum_{k=1}^{\infty}\|B^{-(k+1)}\|\,\|b_1({\widetilde{F}_\ell}^k(x))\!-\!B\|\,\|c_k(u,v,w)-c_k(u,v,\tilde{w})\|
\le M\|w\!-\!\tilde{w}\|^{p_1}.
\label{se-proof2}
\end{align}

In order to prove (\ref{se3}), note that $\|c_1({\widetilde{F}_\ell}^k(u,v,w))\|\le K(\mu_\ell^+)^k$
by (\ref{estim2}) and that
\begin{align*}
\|c_1({\widetilde{F}_\ell}^k(u,v,w))-c_1({\widetilde{F}_\ell}^k(u,v,\tilde{w}))\|
\le L(\mu_m^+)^{k\zeta}\|w-\tilde{w}\|^{{\zeta}}
\end{align*}
by (\ref{estim16-D}). In view of Lemma \ref{lm-hest}, where choosing
\begin{align*}
\alpha=\zeta,~~~~~\tau_1=\mu_\ell^+, ~~~~~ \tau_2=(\mu_m^+)^\zeta,
~~~~~ \rho=\mu_m^+/\mu_\ell^-.
\end{align*}
and noticing that $\rho\tau_1=\eta<1$, we obtain
\begin{align}
&\sum_{k=1}^{\infty}\|B^{-(k+1)}\|\,\|c_1({\widetilde{F}_\ell}^k(u,v,w))
       -c_1({\widetilde{F}_\ell}^k(u,v,\tilde{w}))\|\,\|C^k\|
\nonumber\\
&\le\frac{1}{\mu_\ell^-}\sum_{k=1}^{\infty}\Big(\frac{\mu_m^+}{\mu_\ell^-}\Big)^k\|c_1({\widetilde{F}_\ell}^k(u,v,w))
       -c_1({\widetilde{F}_\ell}^k(u,v,\tilde{w}))\|
\nonumber\\
&\le M\|w-\tilde{w}\|^{p_2},
\label{se-proof3}
\end{align}
where
\begin{align*}
p_2&:=\min\Big\{\zeta-\varepsilon,
~\frac{\log\mu_\ell^++\log\mu_m^+-\log\mu_\ell^-}{\log\mu_\ell^+-\zeta\log\mu_m^+}\,\zeta\Big\}
\\
   &~=\min\Big\{\zeta-\varepsilon, ~\frac{\log\lambda_\ell^++\log\lambda_m^+-\log\lambda_\ell^-}
     {\log\lambda_\ell^+-\zeta\log\lambda_m^+}\,\zeta-\varepsilon\Big\}>0.
\end{align*}
Finally, putting
\begin{align*}
\beta_{\ell}=\min\{p_1,p_2\} =\min\Big\{\zeta-\varepsilon,
~\frac{\log\lambda_\ell^++\log\lambda_m^+-\log\lambda_\ell^-}
     {\log\lambda_\ell^+-\zeta\log\lambda_m^+}\,\zeta-\varepsilon\Big\}>0,
\end{align*}
we prove (\ref{se1})-(\ref{se3}) by (\ref{se-proof1})-(\ref{se-proof3}) respectively.
The proof of the lemma is completed.
\quad$\Box$


\section{$C^{1,\beta}$ linearization of contractions}
\setcounter{equation}{0}

Now we are ready to give the result of $C^{1,\beta}$ linearization of contractions.

\begin{lm}
Suppose that $F:X\to X$ is a
$C^{1,1}$ mapping fixing $O$ and that {\rm (\ref{spec-decomp})}-{\rm (\ref{spec-order})} hold with inequality
{\rm (\ref{NR-conB})}, i.e.,
\begin{align*}
\lambda_i^+/\lambda_i^-<(\lambda_m^+)^{-1}, ~~~~~ \forall i=1,...,m.
\end{align*}
Then there exists a $C^{1,\beta_1}$ diffeomorphism $\Phi: U\to X$ such that equation
{\rm (\ref{schd-eqn})} holds,
where $\beta_1>0$ is the first element of the sequence $(\beta_\ell)_{\ell=1,...,m}$ defined recursively by $\beta_m:=1$ and
\begin{align*}
\left\{
\begin{array}{ll}
\beta_{\ell}:=\min\Big\{\zeta_\ell-\varepsilon, ~\frac{\log\lambda_\ell^++\log\lambda_m^+-\log\lambda_\ell^-}
     {\log\lambda_\ell^+-\zeta_\ell\log\lambda_m^+}\,\zeta_\ell-\varepsilon\Big\},
\\
\vspace{-0.1cm}
\\
\zeta_\ell:=\min\{\beta_{\ell+1},~\log\lambda_\ell^+/\log\lambda_{\ell+1}^- -1-\varepsilon\},
\end{array}
\right.
~~~ \forall \ell=1,...,m-1,
\end{align*}
for any small given $\varepsilon>0$.
\label{lm-holder-B}
\end{lm}


{\bf Proof}. As mentioned in the
beginning of Section 5, we assume
that
$
\Lambda={\rm diag} (\Lambda_1,...,\Lambda_m).
$
For the first step, let
$$
F_m=F=({\rm diag}
(\Lambda_1,...,\Lambda_{m-1})u+\tilde{f}_{a}(u,v),\Lambda_kv+\tilde{f}_{b}(u,v)),
$$
where $u\in X_1\times\cdots\times X_{m-1}$, $v\in X_{m}$,
$\tilde{f}_{a}:X\to X_1\times\cdots\times X_{m-1}$ and
$\tilde{f}_{b}:X\to X_{m}$. Then we can transform
$F_m$ into
$F_{m-1}$, which is in the form of (\ref{F-induc}) with
$\ell=m-1$. In fact, we have the estimate
\begin{align*}
&\|\Lambda_m^{-(k+1)}D(\pi_v F_m^{k+1})(x)-\Lambda_m^{-k}D(\pi_v
F_m^{k})(x)\|
\nonumber\\
&\le \|\Lambda_m^{-(k+1)}\|\,\|D(\pi_v
F_m^{k+1})(x)-\Lambda_mD(\pi_v F_m^{k})(x)\|
\nonumber\\
&\le
\|\Lambda_m^{-(k+1)}\|\,\|DF_m(F_m^k(x))-\Lambda\|\,\|DF_m^k(x)\|
\nonumber\\
&\le M\Big(\frac{(\mu_m^+)^2}{\mu_m^-}\Big)^k
\end{align*}
by (\ref{fn'<uk}) because $F_m$ is
$C^{1,1}$, and the estimate
\begin{align*}
& \|\{\Lambda_m^{-(k+1)}D(\pi_v F_m^{k+1})(x)-\Lambda_m^{-k}D(\pi_v
F_m^{k})(x)\}
\nonumber\\
 &\hspace{2.5cm}   -\{\Lambda_m^{-(k+1)}D(\pi_v F_m^{k+1})(\tilde{x})-\Lambda_m^{-k}D(\pi_v F_m^{k})(\tilde{x})\}\|
\nonumber\\
& \le \|\Lambda_m^{-(k+1)}\|\,\|\{D(\pi_v
F_m^{k+1})(x)-\Lambda_mD(\pi_v F_m^{k})(x)\}
\nonumber\\
 &\hspace{2.5cm} -\{D(\pi_v F_m^{k+1})(\tilde{x})-\Lambda_mD(\pi_v F_m^{k})(\tilde{x})\}\|
\nonumber\\
& \le
\|\Lambda_m^{-(k+1)}\|\,\|DF_m(F_m^k(x))-DF_m(F_m^k(\tilde{x}))\|\,\|DF_m^k(x)\|
\nonumber\\
&
\hspace{2.5cm}+\|\Lambda_m^{-(k+1)}\|\,\|DF_m(F_m^k(\tilde{x}))-\Lambda\|\,\|DF_m^k(x)-DF_m^k(\tilde{x})\|
\nonumber\\
& \le M\Big(\frac{(\mu_m^+)^2}{\mu_m^-}\Big)^k\|x-\tilde{x}\|
\end{align*}
by (\ref{fn'<uk}), (\ref{xy1}) and (\ref{xany1}) with $\zeta=1$.
Since $(\mu_m^+)^2/\mu_m^-<1$ by (\ref{NR-conB}) with $i=m$, one
concludes that the limit $
\Psi_m:=\lim_{k\to\infty}\Lambda_m^{-k}\pi_vF_m^k $ gives a
$C^{1,1}$ mapping. Then, using similar arguments to the second
paragraph of the proof of Lemma \ref{lm-induc}, $F_m$ can be
transformed into $F_{m-1}$ by a $C^{1,1}$ transformation
$\Phi_m(u,v):=(u, \Psi_m(u,v))$ and therefore $F_{m-1}$ is also
$C^{1,1}$ such that $F_{m-1}(O)=O$ and $DF_{m-1}(O)=\Lambda$.

For an inductive proof, assume that $F$ is partially linearized, i.e., can be transformed into the $C^1$
mapping $F_{\ell}$ given in (\ref{F-induc}) such that
(\ref{fafb})-(\ref{hold-l6}) hold for an $\ell\in\{2,...,m-1\}$.
Note that (\ref{NR-conB}) implies (\ref{NR-conB-lll}). Thus,
by Lemma \ref{lm-induc}, $F_\ell$ can be transformed into
a $C^1$ mapping $F_{\ell-1}$ given in (\ref{Fl-1-induc})
by a $C^{1,\beta_{\ell}}$ diffeomorphism $\Phi_{\ell}$. Let
$$
u:=(u^*,v^*)\in (X_1\times\cdots\times X_{\ell-2})\times X_{\ell-1},
~~~~ w^*:=(v,w)\in X_\ell\times\cdots\times X_{m}
$$
and let
$$
{A^*}:={\rm diag}(\Lambda_1,...,\Lambda_{\ell-2}), ~~~~
{B^*}:=\Lambda_{\ell-1}, ~~~~ {C^*}:={\rm
diag}(\Lambda_{\ell},...,\Lambda_{m})
$$
accordingly. Then $F_{\ell-1}$ can be rewritten as
$$
F_{\ell-1}(x)=(A^*u^*+{g}_a(u^*,v^*,w^*), B^*v^*+{g}_b(u^*,v^*,w^*),
C^*w^*),
$$
where ${g}_a: U\to X_1\times\cdots\times X_{\ell-2}$ and ${g}_a:
U\to X_{\ell-1}$ are both $C^1$ and satisfy (\ref{fafb}) where $f$ is replaced by $g$.
Furthermore, in virtue of
(\ref{hold-l+11})-(\ref{hold-l+13}), one can check that $g_a$ and
$g_b$ satisfy (\ref{hold-l1})-(\ref{hold-l6}) by replacing $f$ and
$\beta_{\ell+1}$ with $g$ and $\beta_{\ell}$ respectively.
Therefore, the above inductive proof shows that we can finally get a $C^{1}$ mapping $F_1$ in the form
\begin{align}
F_1=(\Lambda_1v+\tilde{g}(v,w),{\rm diag}
(\Lambda_2,...,\Lambda_{m})w)
\label{F111},
\end{align}
where $v\in X_1$, $w\in X_2\times\cdots\times X_{m}$ and $\tilde{g}:
U\to X_1$ such that $\tilde{g}(O)=0$, $D\tilde{g}(O)=0$ and
\begin{align*}
& \|\partial_{u}\tilde{g}(x)-\partial_{u}\tilde{g}(\tilde{x})\|\le
M_2\|x-\tilde{x}\|, ~~~~ \|\partial_{w}\tilde g(v,w)-\partial_{w}\tilde
g(\tilde{v},w)\|\le M_2 \|v-\tilde{v}\|,
\\
&\|\partial_{w}\tilde g(v,w)-\partial_{w}\tilde g(v,\tilde{w})\|\le
M_2 \|w-\tilde{w}\|^{\beta_2}.
\end{align*}

At last, we use Lemma \ref{lm-induc} again to transform $F_1$ into
$F_0=\Lambda$ via a $C^{1,\beta_1}$ diffeomorphism $\Phi_1$.
Actually, although $F_1$ given in (\ref{F111}) does not have exactly
the same form as $F_\ell$ given in (\ref{F-induc}), Lemma \ref{lm-induc}
still works since one can observe that the term $Au+ f_{a}(u,v,w)$
is not important to the proof of Lemma \ref{lm-induc}.
Consequently, the proof can be completed by putting
$$
\Phi:=\Phi_1\circ\cdots\circ\Phi_m,
$$
a $C^{1,\beta_1}$ diffeomorphism near $O$. The proof is completed.
\quad$\Box$


\begin{re}
{\rm
For contractions, our Lemma \ref{lm-holder-B} extends the results of \cite{El-JFA01} and
\cite{R-S-JDE04} respectively due to a weaker condition on $\sigma(\Lambda)$ and a larger bound of $\beta$.
More precisely, the band condition {(\ref{NR-conB})} given in Lemma \ref{lm-holder-B} is weaker than the one
(\ref{bond-cond}), which was given in \cite{El-JFA01}. In fact,
{(\ref{NR-conB})} can be deduced from (\ref{bond-cond}) but the spectrum
\begin{align*}
\sigma(\Lambda)=[\frac{1}{16}+\delta,\frac{1}{8}]\cup[\frac{1}{8}+\delta,\frac{1}{4}]\cup[\frac{1}{4}+\delta,\frac{1}{2}],
\end{align*}
where $\delta>0$ is sufficiently small, does not satisfy the band
condition although it satisfies {(\ref{NR-conB})}. Moreover, in
\cite{R-S-JDDE06} an analytic contraction on Hilbert space, which dose not satisfy (\ref{NR-conB}),
is given and proved to be not $C^1$ linearizable. This is a counter example which shows the sharpness
of the condition (\ref{NR-conB}).
}
\label{re-comNRcc}
\end{re}


\section{Main results}
\setcounter{equation}{0}

This section is devoted to the main result of the paper.

\begin{thm}
Suppose that $F:X\to X$ is a
$C^{1,1}$ mapping fixing $O$
and that {\rm (\ref{FFbumpR2})}-{\rm (\ref{mu-lam})} hold with
\begin{align}
\left\{
\begin{array}{ll}
\lambda_u^-/\lambda_s^+>\max\{\lambda_u^+,(\lambda_s^-)^{-1}\},
\\
\vspace{-0.3cm}
\\
\lambda_i^+/\lambda_i^-<(\lambda_d^+)^{-1}, ~~~ \forall i=1,...,d, ~~~
\lambda_j^+/\lambda_j^-<\lambda_{d+1}^-, ~~~ \forall j=d+1,...,m.
\end{array}
\right.
\label{NR}
\end{align}
Then there exists a $C^{1,\beta}$ diffeomorphism $\Phi: U\to X$ such that equation
{\rm (\ref{schd-eqn})} holds, i.e., $F$ can be $C^{1,\beta}$ linearized near $O$. Here
$\beta>0$ is defined by
\begin{align*}
\beta=\min\Big\{\beta_1, ~\frac{\log\lambda_{d}^++\log\lambda_m^+-\log\lambda_{d+1}^-}{\log\lambda_{d}^+-\log\lambda_m^+} -\varepsilon,~
\beta_m, ~\frac{\log\lambda_{d+1}^-+\log\lambda_1^--\log\lambda_{d}^+}{\log\lambda_{d+1}^--\log\lambda_1^-} -\varepsilon\Big\}
\end{align*}
for any small given $\varepsilon>0$, where $\beta_1$ is the first element of the sequence $(\beta_i)_{i=1,...,d}$ defined recursively by $\beta_d=1$ and
\begin{align*}
\left\{
\begin{array}{ll}
\beta_{i}=\min\Big\{\zeta_i-\varepsilon, ~\frac{\log\lambda_i^++\log\lambda_d^+-\log\lambda_i^-}
     {\log\lambda_i^+-\zeta_i\log\lambda_d^+}\,\zeta_i-\varepsilon\Big\},
\\
\vspace{-0.1cm}
\\
\zeta_i=\min\{\beta_{i+1},~\log\lambda_i^+/\log\lambda_{i+1}^- -1-\varepsilon\},
\end{array}
\right.
~~~~~ \forall i=1,...,d-1,
\end{align*}
and $\beta_m$ is the last element of the sequence $(\beta_j)_{j=d+1,...,m}$ defined recursively by $\beta_{d+1}=1$ and
\begin{align*}
\left\{
\begin{array}{ll}
\beta_{j}=\min\Big\{\zeta_j-\varepsilon, ~\frac{\log\lambda_j^-+\log\lambda_{d+1}^--\log\lambda_j^+}
     {\log\lambda_j^--\zeta_j\log\lambda_{d+1}^-}\,\zeta_j-\varepsilon\Big\},
\\
\vspace{-0.1cm}
\\
\zeta_j=\min\{\beta_{j-1},~\log\lambda_j^-/\log\lambda_{j-1}^+ -1-\varepsilon\},
\end{array}
\right.
~~~~~ \forall j=d+2,...,m.
\end{align*}
\label{thm-linB}
\end{thm}

{\bf Proof}.
Since $O$ is a hyperbolic fixed point of $F$, it is well known (see e.g. \cite{HPS-book77})
that there exist a stable manifold ${\cal G}_s$ and an unstable manifold ${\cal G}_u$.
They are graphs of $C^{1,1}$ mappings $g_s:X_-\to X_+$ and $g_u:X_+\to X_-$ respectively such that
$g_s(0)=g_u(0)=0$ and $Dg_s(0)=Dg_u(0)=0$. Then, by two $C^{1,1}$ transformations
\begin{align*}
\Theta_1: (x_-,x_+)\mapsto (x_--g_u(x_+),x_+)~~~ {\rm and} ~~~ \Theta_2: (x_-,x_+)\mapsto (x_-,x_+-g_s(x_-)),
\end{align*}
we can straighten up the stable and unstable manifolds. This enables us to assume that
\begin{align}
{\rm Lip}(g_s)={\rm Lip}(g_u)=0,
\label{invmad-lip}
\end{align}
where Lip means the Lipschitz constant.

In order to construct the invariant foliations, we notice that condition (\ref{NR22}) of Lemma \ref{lm-LPeqn}
is satisfied because of the first inequality of (\ref{NR}). Then, let
$q_0:X\times X_-\to X$ be given in Lemma \ref{lm-1b} which is proved to be $C^{1,\beta_s}$ near $O$.
Here $\beta_s>0$ is given in (\ref{def-bs}). Assume that
\begin{align*}
{\cal M}_s(x):=\{x+q_0(x,y_-):y_-\in X_-\}, ~~~~~ \forall x\in X,
\end{align*}
and that the local $C^{1,\beta_s}$ mapping $h_s:X\times X_-\to X_+$ is defined by
\begin{align*}
h_s(x,y_-)=\pi_+(x+q_0(x,y_-)).
\end{align*}
According to \cite[Theorem 3.1]{ChHaTan-JDE97} and its proof, we see that
\begin{description}
\item (B1)
$\pi_-(x+q_0(x,y_-))=y_-$, which implies that ${\cal M}_s(x)$ is just the graph of $h_s(x,\cdot): X_-\to X_+$;
\item (B2)
$h_s(x,x_-)=x_+$, which implies that ${\cal M}_s(x)$ passes through the point $x$;
\item (B3)
either ${\cal M}_s(x)\cap {\cal M}_s(y)=\emptyset$ or ${\cal M}_s(x)={\cal M}_s(y)$ for all $x,y\in X$,
which implies that
\begin{align*}
{\rm if} ~~ h_s(x,y_-)=y_+ ~~~~~~~ {\rm then} ~~ h_s(x,\cdot)=h_s(y,\cdot);
\end{align*}
\item (B4)
$F({\cal M}_s(x))\subset {\cal M}_s(F(x))$, which implies that
\begin{align*}
F(x+q_0(x,y_-))=F(x)+q_0(F(x),\pi_-F(x+q_0(x,y_-))).
\end{align*}
\end{description}
As shown at the beginning of Section 3, conclusions (B1)-(B4) indicate that the set $\{{\cal M}_s(x):x\in
X\}$ forms a stable foliation of $X$ for $F$. Moreover, the first inequality of (\ref{NR})
also implies that $(\lambda_u^-)^{-1}(\lambda_s^-)^{-1}<(\lambda_s^+)^{-1}$. Then,
the same arguments for the inverse $F^{-1}$ give an unstable
foliation $\{{\cal M}_u(x):x\in X\}$ and its corresponding mapping
$h_u:X\times X_+\to X_-$, which satisfies the analogous properties as
presented in (B1)-(B4) and is locally $C^{1,\beta_u}$ with
\begin{align*}
\beta_u=\frac{\log\lambda_{d+1}^-+\log\lambda_1^--\log\lambda_{d}^+}{\log\lambda_{d+1}^--\log\lambda_1^-} -\varepsilon.
\end{align*}

In view of (\ref{invmad-lip}) and (B1)-(B4), one verifies that
the conditions ${\rm (D1)}$-${\rm (D5)}$ and ${\rm (I1)}$-${\rm (I2)}$ given in \cite{Tan-JDE00} hold with respect to $F$.
Hence, by \cite[Theorem 3.1]{Tan-JDE00}, there are
homeomorphisms $\Psi:X\to X$ and its inverse $\Psi^{-1}:X\to X$, which are both $C^{1,\min\{\beta_s,\beta_u\}}$ near $O$, such that the equality
\begin{align}
\Psi\circ F&=F_-\circ \pi_-\Psi+F_+\circ \pi_+\Psi
\label{docp}
\end{align}
holds, where $F_-: X_-\to X_-$ and $F_+:X_+\to X_+$ are defined by
\begin{align*}
F_-=\pi_-F\circ({\rm id}_-+g_s),~~~~~~~ F_+=\pi_+F\circ({\rm id}_++g_u)
\end{align*}
and id$_-$ and id$_+$ denote the identity mappings of $X_-$ and $X_+$ respectively.

Obviously, $DF_-(O)=\Lambda_-$ and $DF_+(O)=\Lambda_+$, where
$\Lambda_-$ and $\Lambda_+$ are given in (\ref{def-FB}). Recall
the spectra $\sigma(\Lambda_-)=\sigma_-$ and
$\sigma(\Lambda_+)=\sigma_+$, which are given in (\ref{spec-cond}).
Since $F_-$ is a $C^{1,1}$ contraction and $F_+$ is a $C^{1,1}$
expansion,
considering the inverse $F_+^{-1}$ in the case of expansion,
by Lemma \ref{lm-holder-B}
we get immediately the following
result:

\begin{lm}
Suppose that the $C^{1,1}$ mappings $F_-, F_+$ are given above
and that {\rm (\ref{spec-distB})} and the inequalities
\begin{align*}
\lambda_d^+\lambda_i^+<\lambda_i^- ~~~ \forall i=1,...,d, ~~~~~~
\lambda_j^+<\lambda_{d+1}^-\lambda_j^- ~~~ \forall j=d+1,...,m,
\end{align*}
hold.
Then $F_-$ and $F_+$ can be linearized by a $C^{1,\beta_1}$ diffeomorphism
$\psi_-:U_-\to X_-$ and a $C^{1,\beta_m}$ diffeomorphism $\psi_+: U_+\to X_+$respectively,
where $\beta_1,\beta_m>0$ are given in Theorem \ref{thm-linB} and
$U_-\subset X_-$, $U_+\subset X_+$ are neighborhoods of the origin.
\label{lm-lince}
\end{lm}

By Lemma \ref{lm-lince},
\begin{align}
\psi_-\circ F_-=\Lambda_-\circ \psi_-, ~~~~~ \psi_+\circ
F_+=\Lambda_+\circ \psi_+.
\label{two-coneqn}
\end{align}
This enables us to define a local $C^{1,\beta}$ mapping $\Phi$ together with its local $C^{1,\beta}$ inverse $\Phi^{-1}$ by
\begin{align*}
\Phi=(\psi_-\circ \pi_-+\psi_+\circ \pi_+)\circ \Psi,~~~~~
\Phi^{-1}=\Psi^{-1}\circ (\psi_-^{-1}\circ \pi_-+\psi_+^{-1}\circ \pi_+),
\end{align*}
where $\beta:=\min\{\beta_1,\beta_s,\beta_m,\beta_u\}$.
Hence (\ref{docp}) and (\ref{two-coneqn}) yield
\begin{align*}
\Phi\circ F&=(\psi_-\circ \pi_-+\psi_+\circ \pi_+)\circ (F_-\circ \pi_-\Psi+F_+\circ \pi_+\Psi)
\\
&=\psi_-\circ F_-\circ \pi_-\Psi+\psi_+\circ F_+\circ \pi_+\Psi
\\
&= (\Lambda_-\circ \psi_- \circ \pi_-+\Lambda_+\circ \psi_+ \circ \pi_+)\circ \Psi
\\
&= \Lambda\circ(\psi_- \circ \pi_-+\psi_+ \circ \pi_+)\circ \Psi
\\
&= \Lambda\circ \Phi.
\end{align*}
The proof is completed. \qquad$\Box$

\begin{re}
{\rm
Our condition (\ref{NR}) is
weaker than the condition in (\ref{RS-NR}), which was given in \cite{R-S-JDDE04}.
In fact, (\ref{RS-NR}) can be deduced from (\ref{NR}) by putting $d=1$ and $m=2$, whereas the example
\begin{align*}
\sigma(\Lambda):=[1/10, 1/6]\cup[2,3]\cup[9,10]
\end{align*}
satisfies (\ref{NR}) but does not satisfy (\ref{RS-NR}).
}
\label{re-comNR}
\end{re}

Note that the finite-dimensional space $\mathbb{R}^n$ always has
bump functions and that the inequalities given in the second row of (\ref{NR})
hold automatically in $\mathbb{R}^n$.
Applying Theorem \ref{thm-linB} to $X=\mathbb{R}^n$, we have the following corollary:

\begin{cor}
Suppose that $F:\mathbb{R}^n\to \mathbb{R}^n$ is
a $C^{1,1}$ mapping fixing $O$
and that the inequalities {\rm (\ref{def-FB})}, {\rm (\ref{mu-lam})} and the inequality
\begin{align}
\lambda_u^-/\lambda_s^+>\max\{\lambda_u^+,(\lambda_s^-)^{-1}\}
\label{NR11}
\end{align}
hold.
Then there exists a $C^{1,\beta}$ {\rm (}$\beta>0${\rm )} diffeomorphism $\Phi: U\to X$ such that equation
{\rm (\ref{schd-eqn})} holds.
\label{cor-linB}
\end{cor}

\begin{re}
{\rm
Condition (\ref{NR11}) can be satisfied by many examples, for instance, the spectrum
\begin{align*}
\sigma(\Lambda):=\{1/10,1/6\}\cup\{2,5,10\}.
\end{align*}
Corollary \ref{cor-linB} extends Belitskii's result given in \cite{Beli-FAA73}, as mentioned
in the Introduction, because one can only conclude $C^1$ but not $C^{1,\beta}$ linearization by \cite[Theorem 2]{Beli-FAA73}.
}
\label{re-speB}
\end{re}


\section{Remarks on sharpness}
\setcounter{equation}{0}

Finally, we discuss on sharpness of the estimate of $\beta$ presented in Theorem \ref{thm-linB}.
This sharpness was given by \cite{Stowe-JDE86} and \cite{WZWZ-JFA01} for planar hyperbolic mappings
in the case that $0<|\lambda_1|<1<|\lambda_2|$ and the case that $0<|\lambda_1|<|\lambda_2|<1$
respectively with counter examples.
We can use the same counter examples to show the sharpness of the estimate of $\beta$ given
in the following by applying our Theorem \ref{thm-linB} to the case $X=\mathbb{R}^2$.

\begin{cor}
Suppose that $F:\mathbb{R}^2\to \mathbb{R}^2$ is a
$C^{1,1}$ mapping fixing $O$ such that
$
DF(O)=
{\rm diag}(\lambda_1, \lambda_2),
$
where
$
0<|\lambda_1|<1<|\lambda_2|.
$
Then $F$ can be $C^{1,\beta}$ linearized near $O$, where
\begin{align}
\beta:=\min\Big\{\frac{\log |\lambda_1|}{\log |\lambda_1|-\log |\lambda_2|}-\varepsilon,
~~ \frac{\log |\lambda_2|}{\log |\lambda_2|-\log |\lambda_1|}-\varepsilon\Big\}\in (0,1]
\label{yyzs}
\end{align}
with any small given $\varepsilon>0$.
\label{cor-C1-R2}
\end{cor}

We note that $\beta_1$ and $\beta_m$ ($m=2$) given in Theorem \ref{thm-linB} are both equal to $1$ in the present case
and that $\lambda_i^+=\lambda_i^-=|\lambda_i|$ for $i=1,2$.
By Theorem \ref{thm-linB} we get
\begin{align*}
\beta
&=\min\Big\{1, ~\frac{\log\lambda_1^++\log\lambda_2^+-\log\lambda_2^-}{\log\lambda_1^+-\log\lambda_2^+} -\varepsilon,~
~\frac{\log\lambda_2^-+\log\lambda_1^--\log\lambda_1^+}{\log\lambda_2^--\log\lambda_1^-} -\varepsilon\Big\}
\\
&=\min\Big\{\frac{\log |\lambda_1|}{\log |\lambda_1|-\log |\lambda_2|}-\varepsilon,
~~ \frac{\log |\lambda_2|}{\log |\lambda_2|-\log |\lambda_1|}-\varepsilon\Big\},
\end{align*}
which proves this corollary.

We assert that the $\beta$ given in (\ref{yyzs}) by our Theorem~\ref{thm-linB} is sharp.
In fact,
putting $\sigma:=-{\log|\lambda_2|}/{\log|\lambda_1|}>0$, we get from (\ref{yyzs}) that
\begin{align*}
\beta=\left\{
\begin{array}{ll}
(1+{\sigma})^{-1}-\varepsilon, & ~~ {\rm if} ~ 0<\sigma<1,
\\
\vspace{-0.3cm}
\\
{\sigma}(1+{\sigma})^{-1}-\varepsilon, & ~~ {\rm if} ~ \sigma\ge 1,
\end{array}
\right.
\end{align*}
an almost the same $\beta$ as given in \cite[Theorem 2]{Stowe-JDE86}
except for the number $\varepsilon$.
Remark that the number $\varepsilon$,
coming from not only the number $\delta>0$ but also the number $\eta>0$ (as seen in
the equality just below (\ref{D3est})), cannot be avoided
because $\eta$ cannot be removed technically
from the inequality $\|DF^k(x)\|\le (|\lambda_2|+\eta)^k$.
The counter example given in \cite{Stowe-JDE86} shows that
this estimate is sharp (or say `almost sharp' because of the small $\varepsilon$).

Similarly, applying our Theorem \ref{thm-linB} to planar $C^{1,1}$
contractions, we also obtain a $C^{1,\beta}$ linearization with
$\beta\in (0,1]$ such that
\begin{align*}
\beta=\left\{
\begin{array}{lllll}
1, & ~~~{\rm if}~ \log_{|\lambda_2|}|\lambda_1|<2,
\\
1-\varepsilon, & ~~~{\rm if}~ \log_{|\lambda_2|}|\lambda_1|=2,
\\
\log|\lambda_2|/(\log|\lambda_1|-\log|\lambda_2|),
       & ~~~{\rm if}~ \log_{|\lambda_2|}|\lambda_1|> 2
\end{array}
\right.
\end{align*}
for any given small $\varepsilon>0$.
This result is exactly the same as the one given in \cite[Theorem 2]{WZWZ-JFA01} (with $\alpha=1$).
Therefore, the counter example given in \cite{WZWZ-JFA01} shows its sharpness.

\begin{re}
{\rm
Our Corollary~\ref{cor-C1-R2} concludes a $C^{1,\beta}$ linearization for planar $C^{1,1}$ hyperbolic mappings with both contraction and expansion.
Note that \cite[Theorem 2]{Stowe-JDE86} considers $C^2$ mappings in the form
$
F(x):=\Lambda x+o(\|x\|^2).
$
Therefore, our Corollary~\ref{cor-C1-R2}
improves \cite[Theorem 2]{Stowe-JDE86}.
}
\end{re}

The above discussion shows that the estimates of the exponent $\beta$ given by our Theorem~\ref{thm-linB}
are sharp in the 2-dimensional case.
We conjecture that it is also sharp in the general infinite-dimensional case.
One expects to give in a Banach space such a counter example which cannot be
linearized by any transformation with more smoothness than the one given in Theorem \ref{thm-linB},
but it is not an easy work in an infinite-dimensional setting.
As seen in \cite{WZWZ-JFA01},
when we constructed counter examples in $\mathbb{R}^2$,
we needed a $C^{\infty}$ function
$u:\mathbb{R}^2\backslash\{O\}\to \mathbb{R}$ such that
\begin{align*}
u(x_1,x_2)=
\left\{
\begin{array}{ll}
1 &\mbox{ if } x_1\ge |x_2|,
\\
0 &\mbox{ if } x_1\le 0,
\end{array}\right.
\end{align*}
$\partial {u}/\partial x_1\ge 0$ and $\|{u}^{(r)}(x)\|\le K\|x\|^{-r}$ for $r=1,2$. This function
can be defined by
\begin{align*}
u(x):=\left\{
\begin{array}{ll}
\int_{-\infty}^{x_1/|x_2|}q(t)dt/\int_{-\infty}^{\infty}q(t)dt, &~~ x_2\ne 0,
\\
0, &~~ x_1<0, ~ x_2=0,
\\
1, &~~ x_1>0, ~ x_2=0,
\end{array}
\right.
\end{align*}
where
\begin{align*}
q(x):=\left\{
\begin{array}{ll}
e^{\frac{1}{t(t-1)}}, &~~ 0<t<1,
\\
0, &~~ \mbox {other}.
\end{array}
\right.
\end{align*}
However, we hardly realize this idea in a general case in Banach spaces.

{\footnotesize

\end{document}